\documentclass[reqno]{amsart}
\usepackage{amsmath}
\usepackage{amsfonts}
\usepackage{amssymb}
\usepackage{amsbsy}
\usepackage{graphicx,color}
\usepackage{algorithm}
\usepackage{algpseudocode}
\usepackage{comment}

\usepackage{mathtools}
\usepackage{bm}
\usepackage{xcolor}
\def\btheta{{\boldsymbol \theta}}

\usepackage[hidelinks]{hyperref}
\hypersetup{
	colorlinks,
	citecolor=red!70!black,
	filecolor=black,
	linkcolor=blue!70!black,
	urlcolor=black
}

\calclayout

\newtheorem{theorem}{Theorem}

\theoremstyle{remark}
\newtheorem{remark}{Remark}
\newtheorem{example}{Example}
\theoremstyle{definition}
\newtheorem{definition}{Definition}
\def\bxi{{\boldsymbol \xi}}
\def\bphi{{\boldsymbol \phi}}
\def\bpsi{{\boldsymbol\psi}}

\def\cd{\cdot}

\def\e{{\mathbf e}}
\def\f{{\mathbf f}}

\def\F{{\mathbf F}}
\def\g{{\mathbf g}}
\def\G{{\mathbf G}}

\def\H{{\mathbf H}}
\def\r{{\mathbf r}}

\def\R{{\mathbf R}}

\def\N{{\mathbf N}}
\def\u{{\mathbf u}}
\def\l{{\mathbf l}}
\def\U{{\mathbf U}}
\def\V{{\mathbf V}}
\def\v{{\mathbf v}}
\def\w{{\mathbf w}}

\def\x{{\mathbf x}}

\def\y{{\mathbf y}}

\def\0{{\mathbf 0}}

\def\eps{{\epsilon}}

\def\gl{{\lambda}}

\def\bPhi{{\mathbf \Phi}}
\def\bphi{{\mathbf \phi}}
\def\bPsi{{\mathbf \Psi}}
\def\bpsi{{\mathbf \psi}}
\def\ba{{\bf a}}
\def\bb{{\bf b}}
\def\bxi{{\mathbf \xi}}

\def\cI{{\mathcal I}} \def\cJ{{\mathcal J}} \def\cK{{\mathcal K}}
 \def\cU{{\mathcal U}}  
 \def\cL{{\mathcal L}} \def\cV{{\mathcal V}} 
 
\def\fr{\frac}

\def\re{I\!\!R}
\def\nt{\notag}

\def\diag{\text{diag}}

\def\bxi{{\boldsymbol \xi}}
\def\bphi{{\boldsymbol \phi}}
\def\bpsi{{\boldsymbol\psi}}

\def\cd{\cdot}

\def\e{{\mathbf e}}
\def\f{{\mathbf f}}
\def\F{{\mathbf F}}
\def\g{{\mathbf g}}
\def\G{{\mathbf G}}
\def\H{{\mathbf H}}
\def\r{{\mathbf r}}

\def\R{{\mathbf R}}
\def\N{{\mathbf N}}
\def\u{{\mathbf u}}
\def\l{{\mathbf l}}
\def\U{{\mathbf U}}
\def\V{{\mathbf V}}
\def\v{{\mathbf v}}
\def\w{{\mathbf w}}

\def\x{{\mathbf x}}

\def\y{{\mathbf y}}

\def\0{{\mathbf 0}}
\def\bxi{{\boldsymbol \xi}} \def\bphi{{\boldsymbol \phi}} \def\bpsi{{\boldsymbol\psi}}
\def\eps{{\epsilon}}

\def\gl{{\lambda}}
\def\bPhi{{\mathbf \Phi}}

\def\bPsi{{\mathbf \Psi}}

\def\ba{{\bf a}}
\def\bb{{\bf b}}

\def\bxi{{\boldsymbol \xi}} \def\bphi{{\boldsymbol \phi}} \def\bpsi{{\boldsymbol\psi}}

\def\cI{{\mathcal I}} \def\cJ{{\mathcal J}} \def\cK{{\mathcal K}}
 \def\cU{{\mathcal U}}  
 \def\cL{{\mathcal L}} \def\cV{{\mathcal V}}

\def\fr{\frac}

\def\re{I\!\!R}
\def\nt{\notag}

\def\diag{\text{diag}}

\begin{document}

\title[Asymmetric exact controllability of coupled elastic strings and springs]{
Asymmetric exact controllability for networks of spatial elastic strings, springs and masses}
\author{G\"unter  Leugering} 
\address{Department Mathematik, Friedrich-Alexander University Erlangen-Nuremberg \\ 91058  Erlangen, Germany}
\email{guenter.leugering@fau.de}

\author{Charlotte Rodriguez}
\address{Department Mathematik, Friedrich-Alexander University Erlangen-Nuremberg \\ 91058  Erlangen, Germany}
\email{charlotte.rodriguez@fau.de}

 \author{Yue Wang $^*$} \thanks{$^*$ Corrseponding author, Fudan University, China (yuewang@fudan.edu.cn); FAU Erlangen-Nürnberg, Germany (yue.wang@fau.de).}
 \address{Department of Data Science,Friedrich-Alexander University Erlangen-Nuremberg \\ 91058  Erlangen, Germany}
 \email{yue.wang@fau.de}

\thanks{ This research was funded by 
DFG WA5144/1-1 (no.504042427) and STCSM no.25ZR1404013.}

\begin{abstract}
We consider networks of elastic strings with end masses, where the coupling is modeled via elastic springs. The model is representative of a network of nonlinear strings, where the strings are coupled to elastic bodies. The coupled system converges to the classical string network model with Kirchhoff and continuity transmission conditions as the spring stiffness terms approach infinity and the masses at the nodes vanish. Due to the presence of point masses at the nodes, the boundary conditions become dynamic, and consequently, the corresponding first-order system of quasilinear balance laws exhibits nonlocal boundary conditions. We demonstrate well-posedness in the sense of semi-global classical solutions \cite{li} (i.e., for arbitrarily large time intervals provided that the initial and boundary data are small enough) and observe extra regularity at the masses as in \cite{WangLeugeringLi2017,WangLeugeringLi2019}. We prove local and global-local exact boundary controllability of a star-like network when control is active at the endpoints of the string-spring-mass system, except for one clamped end. In this case, at multiple nodes, a complex smoothing pattern appears, leading to asymmetric control spaces when the springs and masses are present. Furthermore, the rank of the Laplacian matrix at the junction is crucial for the controllability property, particularly in models containing wave equations with degeneration at dynamic boundaries, which can be interpreted as damage in mechanical vibration systems where parts of the springs are missing. 


\end{abstract}

\maketitle


\section{Introduction}
We consider networks of elastic strings in $\re^3$ connected via elastic springs. The elastic springs are representative of elastic Hookean bodies to which the strings are attached. The Hookean body is then given by a network of springs itself. The nonlinear strings are modeled as parametrized curves and, while non-linear equations for the position vectors along the curve are derived, the main emphasis is on networks close to an equilibrium configuration in which all the strings are {\em stretched}, with displacements which necessarily also are vectors in $\re^3$. 

The interest in such kind of elastic coupling of strings is motivated by recent work on mechanical structures that exhibit degeneration at some point inside or at serial of multiple joints. See e.g. \cite{AlabauCannarsaLeugering2017,KogutKupenkoLeugering2022a,KogutKupenkoLeugering2022b}, where degeneration at a given boundary or internal point is considered. Instead, in the current work, we consider a more particular kind of potential damage, as networks of elastic springs between the structural elements, here nonlinear elastic strings, allow for a more subtle way of describing defects.  

In \cite{s:b}, local well-posedness and controllability results were proved for a single non-linear string governed by the non-linear system.  This was done using results on quasilinear hyperbolic systems. 

In \cite{gl}, Li and Gu proved well-posedness and exact controllability for a model of a planar tree network of strings governed by quasilinear wave equations governing {\em scalar} (necessarily transverse) displacements of the strings.  Analogous results in $\re^3$, where the displacements of the strings from ``stretched'' equilibrium configurations are necessarily in three-space, were provided by Leugering and Schmidt \cite{LeugeringSchmidt2011}. For networks of nonlinear Timoshenko beams see analogous results by Gu, Leugering and Li \cite{GuLeugeringLi2017}.

We stay with the concept of semi-global classical solutions developed by Li ~\cite{li}, where semi-global in time means that for any time $T > 0$, and for small enough initial and boundary data, a unique solution exists at least until time $T > 0$. For a similar model of nonlinear strings see \cite{LiSerre,LiPeng,Li}.
In this article we intend to embark on networks of coupled strings as described above. The relaxation of the node conditions appears to be new and is, in fact, physically reasonable. Indeed, exact mathematical models of junctions can be seen as an abstraction from the model provided here.

For the controllability of a system of coupled quasilinear wave equations (of size $n$) with general nonlinear dynamic boundary conditions has been established when the system has full controls at every node ($2n$ controls), or at all extremes of the networks ($n$ controls) in a series of works \cite{WangLeugeringLi2017, WangLeugeringLi2019,LeugeringLiWang2019, WangLeugeringLi2020, LiWang2018}. In these works, a general framework for handling higher-order differential boundary operators was developed, and local exact boundary controllability for various mass-viscoelastic spring-string networks with real-world applications was demonstrated using a modular constructive method. Recently, an interesting smoothing pattern for the solution and the spaces of controllable initial data for two coupled linear elastic strings was shown in \cite{LeugeringMicuRoventaWang2021}. It was found that singularities in waves are "smoothed by three orders" as they cross a point mass when one $L^2$ control acts on the extremity of the first string. In this article, we aim to establish the controllability for networked quasilinear string-spring-mass systems using $n$ or $n-1$ boundary controls at the extremities, and to obtain the asymmetric controllable spaces for semi-global classical solutions through characteristic analysis and constructive methods.

\subsection{Organization of this paper and main results}
This article is organized as follows. In Section \ref{sec:modeling}, we derive the mathematical model for our network of strings and elastic bodies, and in particular the transmission conditions at junctions. This leads to the System \eqref{full-nonlinear-sys} presented below. In Section \ref{sec:wellposedness}, after having presented an star-shaped network, we prove that there exists a unique semi-global solution to \eqref{full-nonlinear-sys}, yielding our first main result Theorem \ref{existence-semi-global}. This result will notably be of use in Section \ref{sec:control}, where is the specific case of a star-like network of $n$ strings controlled at $n-1$ extremities, say all endpoints except for one clamped end, we prove local exact controllability of the network, our second main result Theorem \ref{local-exact-star}. We conclude with Theorem \ref{global-local}, proving a ``global-local'' exact controllability result for star-like networks, in the sense that we may steer the solution to ``far away'' states if the starting and arriving states lay close to some stretched equilibrium solutions. At the end of section  \ref{sec:control}, we consider the controllability properties for the model with damage at junction, which occurs in the form of missing springs at the joints. In two examples, we show that the rank condition is a sufficient but not necessary condition for the controllability property.

\subsection{Notation}
We end this introduction with some comments on notation. Vectors, or vector valued functions, will be indicated in boldface. For a vector $\v$ in a Euclidean space we let $|\v|$ denote the Euclidean length. We denote the Fr\'echet derivative of a function $f$ with respect to a scalar or vector argument $\xi$ by $D_\xi$.  We shall also often write $f_s$ for the partial derivative of $f$ with respect to a scalar variable $s$ and $Df$  for the Fr\'{e}chet derivative of $f$ with respect to its complete argument.
We let $C^n([0,L];W)$ denote the space of $n$ times continuously differentiable functions $\f(x)$ from the interval $[0,L]$ to an open subset W of a Euclidean space, with corresponding norms 
\begin{align*}
\|\f\|_0 = \|\f\| := \sup_{x\in[0,L]} |\f(x)|,\ \ \ \|\f\|_n:=\max\{\|\f\|,\|\f_x\|,\dots,
\|D_x^n\f\|  \}.
\end{align*}
Similar notation will be used for spaces of functions of variables
$(x,t)\in [0,L]\times[0,T]$ or of $t\in[0,T]$ and the associated
norms.


\section{Modeling of networks of nonlinear elastic strings coupled by springs}
\label{sec:modeling}

\begin{figure}
\includegraphics[scale=0.7]{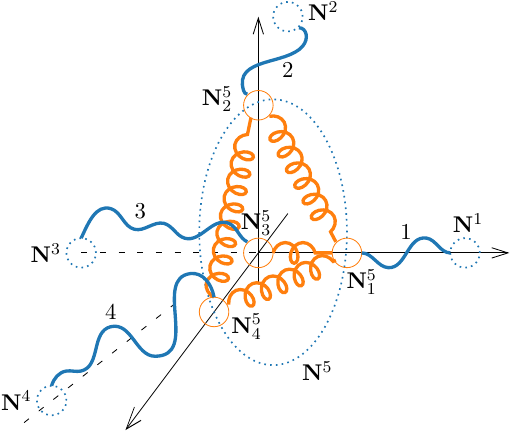}
\hspace{0.75cm} 
\includegraphics[scale=0.75]{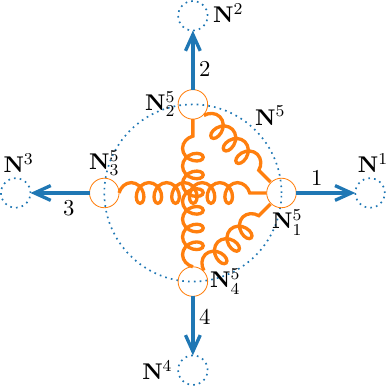}
\caption{A star-like network consisting of four strings represented in $\re^3$ (left) and as a diagram (right).}
\label{fig:3d_and_notation}
\end{figure}

In this section we  describe a nonlinear model for networks of elastic strings and springs, and refer the reader to Figure \ref{fig:3d_and_notation} for visualization.  We suppose that there are $n$ strings indexed by $i \in \cI  = \left\{1,\cdots,n\right\}$. We let the $i-$th string be parameterized by its rest arc length $x$ with $x\in[0,L_i]$, $L_i$ of course being the natural length of that string. The position \textcolor{black}{in $\re^3$} at time $t\in [0, T]$ of the point corresponding to the parameter $x$ will be denoted by the vector $\R^i(x,t)$.  We shall write $\R$ to denote $\{ \R^i\}_{i\in\cI}$. The positions of the endpoints, which we refer to as nodes, are given by functions $\N^j \colon [0, T] \rightarrow \re^3$ with $j \in \cJ  = \left\{1,\cdots,m\right\}$. At time $T>0$, this network is thus described by a graph $\mathcal{G}(t)$ with nodes $\mathcal{N}(t) = \{\mathbf{N}^j(t) \colon j \in \mathcal{J}\}$.

At {\em multiple nodes}, where several strings meet, there is a common location $\N^j(t)$, in fact, there is a common region itself represented by a network of springs, as introduced below. {\em Simple nodes} are those corresponding to the endpoints of only one string. We let 
\begin{align*}
\cI^j = \{i \in \cI \colon \N^j \mbox{ is an end point of the i-th string} \},
\end{align*}
$\cJ^M$ be the subset of $\cJ$ corresponding to
multiple nodes, while $\cJ^S$ contains the indices of simple nodes. We assume that there are simple nodes so that $\cJ^S$ is not empty. \textcolor{black}{We also partition the set of simple nodes as $\cJ^S = \cJ^S_D \cup \cJ^S_N$, the subset $\cJ^S_D$ corresponding to the nodes at which Dirichlet boundary conditions are applied, while $\cJ^S_N$ corresponds to nodes at which Neumann boundary conditions are applied, with boundary data denoted $\mathbf{U}^j(t)$.} \textcolor{black}{The set of indices of nodes incident with the string $i$ is denoted $\cJ^i = \{j_1, j_2 \colon i \in \cI^{j_1} \cap \cI^{j_2}\}$.}
Note that for $j \in \cJ^S$ we let $\cI^j =\{i_j\}$. For $i \in \cI^j$ we let
\begin{align*}
x_{ij} = \begin{cases}
0 &\text{if } \R^i(0,t) = \N^j(t)\\
L_i &\text{if } \R^i(L_i,t) = \N^j(t).
\end{cases}
\end{align*}
For purposes of integration by parts we
also introduce $\eps_{ij}$ to equal $1$ or $-1$ depending on whether $x_{ij}$ is equal to $L_i$ or 0, respectively.  Then, $\eps_{ij}\R^i_x$ is the outward pointing derivative at the boundary point $x_{ij}$ of the interval $[0,L_i]$. 

Instead of coupling the strings at the node $\N^j(t)$, which corresponds to the classical situation (see \cite{LeugeringSchmidt2011}), we expand this node into a graph \textcolor{black}{$\mathcal{G}^j(t)$} itself, the nodes of which, denoted by $\mathcal{N}^j(t):=\{ \N^j_i(t) \colon i=1,\dots, d_j\}$, play the role of $\N^j(t)$. That is to say, if the ``multiple node'' in the classical sense has an edge degree $d_j$ then the graph inserted has $d_j$ nodes and each node is precisely connected to one of the strings. The ``local'' spring network is described by \textcolor{black}{$A^j=(a_{\alpha \beta}^{j})_{\alpha, \beta =1,\dots, d_j}$},  the local adjacency matrix of the inserted graph $\mathcal{G}^{j}$, i.e. $a_{\alpha \beta}^{j}=1$ if there is a connection between the node $\N^{j}_\alpha$ and $\N_\beta^{j}$, otherwise $a^{j}_{\alpha \beta}=0$. We now assume that the connections within the inserted graph $\mathcal{G}^{j}$ are realized by elastic springs with stiffness $\kappa_{j}$. This corresponds to an isotropic elastic body to which the strings are connected.  Clearly, we could have springs with different stiffness and we could consider more general graph insertions such that the graph $\mathcal{G}^{j}$ has itself internal multiple nodes. For the sake of clarity of the exposition, we confine ourselves to the situation above. In this paper, we consider Hookean elastic springs governed by Hooke’s law with given stiffness parameter, even though the fully nonlinear spring can be handled similarly.  We assume that the endpoint of the string $i$ carries a mass $m_{i}^{j}$ at the node $j$.

\begin{figure}
\centering
\includegraphics[scale=0.75]{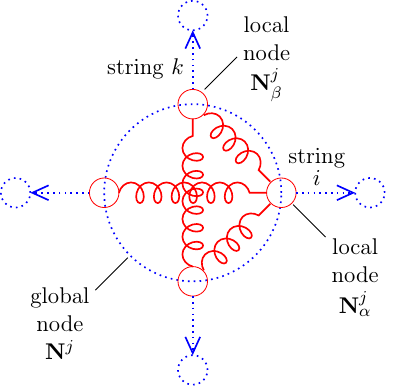}
\caption{A consistent indices notation for the global string graph $\mathcal{G}$ and the local spring graph $\mathcal{G}^j$ -- Visualization with a star-like network.}
\label{fig:notation_indices}
\end{figure}

\begin{remark}
\color{black} To ease the reading, we use the following notation throughout this article. The letters $i, k$ are used for indexes in $\mathcal{I}$ referring to strings, and $j \in \mathcal{J}$ is used for indexes of nodes of the ``global'' graph $\mathcal{G}$. On the other hand, for any given $j \in \mathcal{J}$, we use the letters $\alpha, \beta$ for indexes in $\{1, \ldots, d_j\}$ of nodes of the ``local'' spring graph $\mathcal{G}^j$.
In addition, for the adjacency matrix $A^j$, we will make an abuse of notation, writing $a_{ik}^j$ where $i, k \in \mathcal{I}^j$ to refer to the quantity $a_{\alpha, \beta}^j$ where $\alpha, \beta \in \{1, \ldots, d_j\}$ are the indexes of the nodes of the spring graph which are incident with the strings $i,k$, respectively. We refer to Fig. \ref{fig:notation_indices} for visualization.
\end{remark}

Let $\rho_i$  be the constant density of the corresponding string. Then, the kinetic energy of a
single string, labelled by $i$, at time t is given by
\begin{equation*} 
\cK^i(\R^i(\cd,t)) := \tfrac{1}{2} \int^{L_{i}}_0 \rho_i |\R_t^i (x,
t) | ^2 \, dx+ \frac{1}{2} \textcolor{black}{\sum_{j \in \cJ^i\cap \cJ^M} m_{i}^{j}}|\R^{i}_t(x_{ij},t)|^{2}.
\end{equation*}
We shall assume that the potential energy of the same string is of
the form
\begin{align*}
\cV^i(\R^i(\cd,t))&:=\int^{L_{i}}_{0}\left[ V^i(|\R _x^i (x,
t)|)+\rho_i g\R^i(x,t)\cd\e\right]\,dx\\
&\quad \quad + \frac{1}{2}\textcolor{black}{\sum_{j\in\mathcal{J}^i \cap \cJ^M} \sum\limits_{k \in \mathcal{I}^j, \, k>i}} \kappa_{j}a_{ik}^{j}|\R^{i}(x_{ij},t)-\R^{k}(x_{kj},t)|^{2}\nt,
\end{align*}
where
\begin{itemize}
\item  $V^i(s)$  is a twice continuously differentiable, convex
real valued function defined on an  open subinterval $I^i= (a^i,b^i)$ of the positive real axis, with $ a^i<1<b^i$, satisfying $V^i_{ss}(s) > 0$ and $V^i(1) = V^i_s(1) = 0$; 
\item $\e$ is the vertical unit vector and $g$ is the gravitational constant.
\end{itemize}
As for the total kinetic energy and total potential energy, we define
\begin{equation*}
\cK(\R(\cd,t)) := \sum\limits_{i\in \mathcal{I}}\cK^i(\R^i)
\end{equation*}
and
\begin{equation*}
\cV(\R(\cd,t)):=\sum\limits_{i\in
\mathcal{I}}\cV^i(\R^i),
\end{equation*}
respectively.

\begin{remark}
\textcolor{black}{The first term in the elastic energy corresponds to the Hookean law, accounting for extension or compression — as measured by $|R^i_x(x)|$ — together with the contribution from gravity. While the second term corresponds to potential energy provided by the springs $\kappa_ja^{j}_{ik}|\R^{i}(x_{ij})- \R^{k}(x_{kj})|$ according to Hooke's law.} The generic hypotheses on $V^i$ are quite broad and in the nature of minimum physically plausible assumptions. At points where $|\R_x^i(x)| =1$ there is neither compression nor extension and hence there should be no contribution to the first term in the potential energy while at points where there is compression or extension there is a positive contribution. This leads to the condition that $V^i(s)$ takes minimum value  $0$ at $1$. Certainly the potential energy should increase with increasing extension ($|\R_x^i(x)| > 1$ and increasing) or compression ($|\R_x^i(x)| < 1$ and decreasing). The convexity assumption seems appropriate at least over a small interval of values of $|\R_x^i(x)|$.  In \cite{lls},\cite{ls} and \cite{s:a} these functions took the form   $V^i(s) = h^i(s-1)^2/2$ for $s \in I=(0,\infty)$.
\end{remark}

\begin{definition}
The string is respectively {\em stretched}, {\em limp} or {\em compressed at} $x$ depending on whether $|\R^i_x(x)| > 1,\ = 1$ or $< 1$.  We say that the {\em string is stretched} if $|\R^i_x(x)|>1$ for all $x\in [0,L_i]$ and compressed if $|\R^i_x(x)|<1$ for all $x\in [0,L_i]$.
\end{definition}
Physically, one expects significantly different phenomena in these three situations and transitions between them are certainly complicated.

We now  apply Hamilton's principle to the Lagrangian functional defined for fixed $T > 0$ by
\begin{align*}
&{\cL} (\R)  := \int^T_0 \textcolor{black}{\sum\limits_{i\in
\mathcal{I}} \bigg[ } \int^{L_{i}}_0\left[\tfrac{1}{2}\rho_i |\R_t^i (x,t) |^2dx- V^i(|\R^i_x (x, t) | )-\rho_i g\R^i(x,t)\cd\e  \right] \, dx \nt \\ 
& + \frac{1}{2} \textcolor{black}{\sum_{j\in\cJ^i \cap \cJ^M} m_{i}^{j}} |\R^{i}_t(x_{ij},t)|^{2} - \frac{1}{2}\textcolor{black}{\sum_{j\in\cJ^i \cap \cJ^M}\sum\limits_{k \in \mathcal{I}^j, \, k>i}} \kappa_{j}a_{ik}^{j}|\R^{i}(x_{ij},t)-\R^{k}(x_{kj},t)|^{2}\bigg] \, dt. \nt
\end{align*}
This requires $\cL$ to be stationary at $\R$, where
the domain of $\cL$ consists of those $\R$ whose component functions $\R^i$ are in \textcolor{black}{$C^2([0,L_i]\times[0,T];\re^3)$, and} satisfy given initial conditions
\begin{equation}\label{sys-ic}
\R^i(x,0) = \R^{0,i}(x)\ \  \text{and}\ \ \R^i_t(x,0) = \R^{1,i}(x), \ \ \text{ for } x \in [0, L_i],
\end{equation}
as well as {\em prescribed Dirichlet boundary conditions} at simple nodes
\begin{equation}\label{sys-dirichlet}
 \R^{i_j}(x_{i_jj},t) = \U^j(t)\ \ \text{ for } j \in \cJ^S_D\ \  \text{ and  for}\  t\in [0,T].
\end{equation}
We note that, in contrast to the more classical modeling, we do \textbf{\textit{not}} require the {\em continuity condition} at multiple nodes
\begin{equation}\label{nodalconti}
\R^{i_1}(x_{i_1j},t)=\R^{i_2}(x_{i_2j},t)\ \ \text{ for } j\in \cJ^M,\ \  i_1,i_2 \in \cI^j\ \text{ and for }t\in [0,T],
\end{equation}
due to the fact that the end of the strings are coupled via the springs to the other ends of the strings that are adjacent in $\mathcal{G}^{j}$.
Now, one can consider a perturbation $\R^i(x,t) + \gl\r^i(x,t)$ for each $i \in \cI$, where $\r^i(x,t)\in C^2([0,L_i]\times[0,T];\re^3)$ satisfy homogeneous initial conditions and vanish at the simple nodes.  

A necessary condition for $\cL$  to be stationary at $\R$ is that $\left. D_\gl\cL(\R+\gl\r)\right|_{\gl=0} = 0$, where some computations yield
\begin{align*}
&\left. D_\gl\cL(\R+\gl\r)\right|_{\gl=0}=\\
& \int^T_0 \Bigg\{ \ \
\sum\limits_{i\in \mathcal{I}} \ \, \Bigg[ \int^{L_i}_0 \big( \rho_i \R^i_t \cd \r_t^i-V^i_s(|\R^i_x|)|\fr{\R^i_x}{|\R^i_x|}\cd\r
^i_x-\rho_i g\r^i\cd\e \big)(x,t)dx \Bigg]\\ 
& \qquad +\sum\limits_{j\in \cJ^M}\Bigg[\sum\limits_{i \in \mathcal{I}^j} m_{i}^{j}\R^{i}_{t}(x_{ij},t)\cdot \r^{i}_t(x_{ij},t)\\
& \qquad \qquad  -\kappa_j\sum_{i=1}^{d_j} \sum\limits_{k \in \mathcal{I}^j, \, k>i} a_{ik}^{j}\left(\R^{i}(x_{ij},t)-\R^{k}(x_{kj},t)\right)\cd \left(\r^{i}(x_{ij},t)-\r^{k}(x_{kj},t))\right)\Bigg] \Bigg\}dt.
\end{align*}
First we choose perturbations with compact support in $(0,L_i)\times(0,T)$.  One can then integrate by parts and obtain in the usual way
 the following non-linear partial differential
equation:
\begin{equation}\label{sys-eq}
\rho_i \R^i_{tt}(x,t) = [\G^i(\R^i_x(x,t))]_x  -\rho_i g \e, \ \ \text{
for each } i \in \cI
\end{equation}
with $\G^i:\re^3\mapsto \re^3$ defined by
\begin{equation*}
\G^i(\v) := V^i_s(|\v|)\fr{\v}{|\v|}.
\end{equation*}
Next, for $j\in \cJ^M$ we choose perturbations with $\r^i = \0$ for $i \notin \cI^j$ and with support in a small neighbourhood of $x_{ij}$ for $i \in \cI^j$.  As there is no continuity condition across the joints, we are led to the multiple node condition
\begin{align}\label{sys-nodal}
&\eps_{ij}\G^i(\R^i_x(x_{ij},t))+m^{j}_{i}\R^{i}_{tt}(x_{ij},t)\\ 
&+ \kappa_j\bigg[ \bigg(  \sum\limits_{k \in \mathcal{I}^j} a^{j}_{ik}\bigg)\R^{i}(x_{ij},t)- { \sum\limits_{k \in \mathcal{I}^j}} a^{j}_{ik}\R^{k}(x_{kj},t)\bigg]=0\ \ \text{ for each } j \in \cJ^M, \ i\in \cI^j.\nonumber
\end{align}

\begin{remark}
We notice that if we add \eqref{sys-nodal} over all string indices $i\in \cI^j$, we obtain
\begin{align*}
\sum\limits_{i\in \cI^{j}}\eps_{ij}\G^i(\R^i_x(x_{ij},t))+m^{j}_{i}\R^{i}(x_{ij},t)=0,
\end{align*}
which is the nodal condition for Kirchhoff-type coupling with mass. And if we then neglect the mass, the classical Kirchhoff condition obtains. Moreover, if we let $\kappa^{j}$ tend to $\infty$, then the classical continuity conditions across the joint appear and we are left with network discussed in \cite{LeugeringSchmidt2011}.
\end{remark}

\begin{remark}
In the formulation of system, we mainly impose Dirichlet boundary conditions \eqref{sys-dirichlet} at the simple nodes $j \in \cJ^S_D$. However, depending on the physical modeling context, it is also possible to replace Dirichlet conditions at some simple nodes by Neumann-type boundary conditions. In such cases, for $j \in \cJ^S_N$ the boundary condition reads
\[
   \eps_{i_jj} G^{i_j}_j\!\left(\R^{i_j}_x(x_{i_jj}, t)\right) = \U_j(t), 
   \qquad t \in (0,T), \; i_j \in I_j.
\]
This reflects situations where forces or fluxes rather than displacements are controlled at certain ends of the network.
\end{remark}

We collect the equations and nodal conditions and write down the entire system as follows:
\begin{equation}\label{full-nonlinear-sys}
\tag{E}
\begin{dcases}
\rho_i \R^i_{tt}(x,t) = [\G^i(\R^i_x(x,t))]_x  -\rho_i g \e,
&\text{in }(0, L_i)\times(0,T), \ \  i \in \cI\nt \\
\R^{i_j}(x_{i_jj},t) = \U^j(t),
&t\in (0,T), \, j \in \cJ^S_D \nt \\
\eps_{i_jj}\G^{i_j}(\R^{i_j}_x(x_{i_jj},t))= \U^{j}(t),
&t\in (0,T), \,  j \in \cJ^S_N\nt \\
\eps_{ij}\G^i(\R^i_x(x_{ij},t))+m^{j}_{i}\R^{i}_{tt}(x_{ij},t)
\\ 
\, + \kappa_j\bigg[ \bigg( \sum\limits_{k \in \mathcal{I}^j} a^{j}_{ik}\bigg)\R^{i}(x_{ij},t)- \sum\limits_{k \in \mathcal{I}^j} a^{j}_{ik}\R^{k}(x_{kj},t)\bigg]=0,
&t\in (0,T), \, j \in \cJ^M\nt, \, i \in \cI^j\\
(\R^i, \R^i_t)(x,0) = (\R^{0,i},\R^{1,i})(x),
&x\in (0,L_i), \, i \in \cI .\nt
\end{dcases}
\end{equation}
We consider stretched equilibria $\R^{e,i}(x), \; i\in\mathcal{I}$, satisfying
\begin{equation}\label{equilibrium-nonlinear-sys}
\tag{Eq}
\begin{dcases}
[\G^i(\R^{e,i}_x(x))]_x  =\rho_i g \e,
& x\in(0, L_i), \, i \in \cI \\
\R^{e,i_j}(x_{i_jj}) = \U^j,
&j \in \cJ^S_D \nt \\
\eps_{i_jj}\G^{i_j}(\R^{e,i_j}_x(x_{i_jj}))= \U^{j},
&j \in \cJ^S_N \\
\eps_{ij}\G^i(\R^{e,i}_x(x_{ij}))\\
\, + \kappa_j\bigg[ \bigg( \sum\limits_{k \in \mathcal{I}^j} a^{j}_{ik}\bigg)\R^{e,i}(x_{ij}) - \sum\limits_{k \in \mathcal{I}^j} a^{j}_{ik}\R^{e,k}(x_{kj})\bigg] = 0,
&j \in \cJ^M\nt, \, i\in \cI^j.
\end{dcases}
\end{equation}
\begin{remark}
In the sequel we concentrate on Dirichlet controls at simple nodes, i.e. $\cJ^S_N=\emptyset$, and $\mathcal{J}^S = \mathcal{J}_D^S$.
\end{remark}
\begin{remark}
In this article, we do not attempt to analyze the general existence of equilibria for system \eqref{full-nonlinear-sys}, which can be quite complicated (see \cite{LeugeringSchmidt2011} for a detailed discussion in the context of nonlinear string networks). Instead, we restrict attention to particular network configurations (such as star-like or chain-like networks considered in Section~4), for which equilibrium states $\R^{e,i}(x), \, i\in\cI$ can be explicitly identified and employed in the controllability analysis. If gravitational effects are neglected, the trivial zero solution 
serves as the equilibrium. 

\end{remark}


\section{Well-posedness of \texorpdfstring{\eqref{full-nonlinear-sys}}{systE}}
\label{sec:wellposedness}

In order to find solutions $\R^i(x,t)$ in
$C^2([0,L_i]\times [0,T];\re^3)$ we need
\begin{equation*}
\begin{aligned}
&\U^j(t) \in C^2([0,T];\re^3),\ \ \text{for }j \in \cJ^S,\\
&\R^{0,i}(x) \in C^2([0,L_i];\re^3),\ \ \text{for }i\in\cI,\\  
&\R^{1,i}(x) \in
C^1([0,L_i];\re^3),\ \ \text{for }i\in\cI,
\end{aligned}
\end{equation*}
and that these functions satisfy the following $C^2$ compatibility conditions:

\begin{eqnarray}\label{compatibility-1}
&&\eps_{ij}\G^i(\R^{0,i}_x(x_{ij}))+\frac{m_i^j}{\rho_i}[\G^i(\R_x^{0,i}(x_{ij}))]_x-g\e\nt
\\  && \quad + \kappa_j\bigg[ \bigg( \sum\limits_{k \in \mathcal{I}^j} a^{j}_{ik}\bigg)\R^{0,i}(x_{ij})-\sum\limits_{k \in \mathcal{I}^j} a^{j}_{ik}\R^{0,k}(x_{kj})\bigg]=0,\ \ \text{ for } j \in \cJ^M\nt, \, i \in \cI^j \\
&& \eps_{ij}\G_\v^i(\R^{0,i}_x(x_{ij}))\R^{1,i}_x(x_{ij}) +\frac{m_i^j}{\rho_i}[\G_\v^i(\R_x^{0,i}(x_{ij}))]_x\R^{1,i}_x(x_{ij})\\  && \quad + \kappa_j\bigg[ \bigg( \sum\limits_{k \in \mathcal{I}^j} a^{j}_{ik}\bigg)\R^{1,i}(x_{ij})-\sum\limits_{k \in \mathcal{I}^j} a^{j}_{ik}\R^{1,k}(x_{kj})\bigg]=0,\ \ \text{ for } j \in \cJ^M\nt, \, i \in \cI^j
\end{eqnarray}
as well as
\begin{equation}\label{compatibility-2}
\begin{aligned}
&\U^j(0)=\R^{0, i_j}(x_{ij}), \ \ \text{ for } j\in\cJ^S,\\
&\U^j_t(0)=\R^{1, i_j}(x_{ij}),\ \ \text{ for } j\in\cJ^S,\\
&\U^j_{tt}(0)=(\rho_{i_j})^{-1}\big[\G^{i_j}\big(\R^{0, i_j}_x \big)\big]_x\left.\right|_{x=x_{ij}} - g\e,\ \
\text{ for } j\in\cJ^S,
\end{aligned}
\end{equation}
where $\cI^j = \{i_j\}$ for $j\in \cJ^S.$ After specifying the
initial data $(\R^{0,i},\R^{1,i})\in C^2([0,L_i];\re^3)\times
C^1([0,L_i];\re^3)$ for all $i$, the
control data $\U^j$ at simple nodes can be uniquely represented by
\begin{align*}
\U^j(t)=\overline{\U}^j(t) +
\R^{i_j,0}(x_{ij})+t\R^{i_j,1}(x_{ij})+t^2/2\,(\rho_{i_j})^{-1}\left[\G^{i_j}(\R^{i_j,0}_x)\right]_x\left.
\right|_{x=x_{ij}},
\end{align*}
where  $\overline{\U}^j(t)$  belongs to
\begin{align*}
C^2_0([0,T];\re^3):= \left\{\overline{\U}\in
C^2([0,T];\re^3) \colon \overline{\U}(0)=\overline{\U}_t(0)=\overline{\U}_{tt}(0)=\0\right\}.
\end{align*}
We shall prove a local existence theorem for \eqref{full-nonlinear-sys} close to a specified stretched equilibrium $\R^e =\{
\R^{e,i}\}_{i\in\cI}$. To this end, we transform the above system into an equivalent first-order system of coupled  quasilinear hyperbolic  equations and then apply the existence theorem for semi-global solutions of Li \cite{Li,WangLeugeringLi2017,WangLeugeringLi2019}.  As an intermediate step
we introduce perturbations away from the given equilibrium by setting
\begin{equation*}
\r^i(x,t) := \R^i(x,t)-\R^{e,i}(x).
\end{equation*}
Noting that the $\{\R^{e,i}\}_{i\in \mathcal{I}}$ do not depend on $t$, the system \eqref{full-nonlinear-sys} is equivalent to
\begin{equation}\label{sys-eqi}\tag{E''}
\begin{dcases}
\rho_i \r^i_{tt}(x,t) = [\G^i(\R^{e,i}_x(x)+\r^i_x(x,t))]_x-\rho_i g\e,
&\text{in }(0, L_i) \times (0, T), \, i\in\cI\\
\r^{i_j}(x_{i_jj},t) = \u^j(t),
&t \in (0, T), \, j\in \cJ^S\\
\eps_{ij}\G^i(\R^{e,i}_x(x_{ij})+\r^{i}_x(x_{ij},t))+m^{j}_{i}\r^{i}_{tt}(x_{ij},t)\\
\  + \kappa_j \bigg[ \bigg(\sum\limits_{k \in \mathcal{I}^j} a^{j}_{ik}\bigg)(\R^{e,i}(x_{ij})+\r^{i}(x_{ij},t)) \\ 
\quad\quad \ - \sum\limits_{k \in \mathcal{I}^j} a^{j}_{ik}(\R^{e,k}(x_{kj})+\r^{k}(x_{kj},t))\bigg]=0,
&t \in (0, T), \, j \in \cJ^M, \, i \in \cI^j\\
(\r^i, \r^i_t)(x,0) =(\r^{0,i}, \r^{1,i})(x),
&x \in (0, L_i), \, i\in\cI,\\
\end{dcases}
\end{equation}
where we have set $\r^{0,i} := \R^{0,i}-\R^{e,i}$ and $\r^{1,i}:=\R^{1,i}$, as well as $\u^j = \U^j - \R^{e, i_j}(x_{i_j j})$.
\begin{remark}\label{vectorial}
Consider a multiple node such that $x_{ij}=0$ for all edges adjacent to $\N^j$. At this node we now have a graph $\mathcal{G}^j$ of springs with adjacency structure given by $A^j=(a_{\alpha, \beta}^j)_{\alpha, \beta = 1, \ldots, d_j}$. We introduce the diagonal matrix \textcolor{black}{(the \emph{degree matrix})}
\begin{align*}
D^j:=\text{diag}\bigg(\bigg(\sum\limits_{\beta=1}^{d_j}a_{\alpha \beta}^j\bigg), \ \alpha=1,2,\dots, d_j\bigg).
\end{align*}
The weighted \emph{discrete Laplacian} (or \emph{Laplacian matrix}) of $\mathcal{G}^j$ is given by ${\bf L}^j:=D^j-A^j$. We also define the diagonal \emph{mass matrix} ${\bf M}^j:=\text{ diag}(m_\alpha^j, \ \alpha=1, \dots, d_j)$ \textcolor{black}{and the variables
\begin{align} \label{eq:notation_rbarj}
\overline{\r}^j := (\r^1, \r^2, \ldots, \r^{d_j})^T, \quad \overline{\R}^{e,j} := (\R^{e,1}, \R^{e,2}, \ldots, \R^{e,d_j})^T,
\end{align}
having values in $\re^{d_j \times 3}$, with a similar definition for $\bar{\r}_x^j$ and $\bar{\r}_{tt}^j$.} 
Moreover, we define the strain matrix
\begin{align*}\color{black}
{\mathbf{K}}^j(\overline{\mathbf{r}}^j_x(0,t)):=\big(\mathbf{G}^1(\R^{e,1}(0)+\r^1_x(0,t)), \, \ldots, \, \mathbf{G}^{d_j}(\R^{e,d_j}(0)+\r^{d_j}_x(0,t)) \big)^T.
\end{align*}
Then, the nodal condition coupling the adjacent strings to the springs can be written in a more concise manner as
\begin{align*}
\color{black}
{\bf M}^j\overline{\r}^j_{tt}(0,t)={\mathbf{K}}^{j}(\overline{\r}^j_x(0,t))-\kappa_j {\bf L}^j(\overline{\mathbf{R}}^{e,j}(0)+\overline{\mathbf{r}}^j(0,t)).
\end{align*}
\end{remark}

\begin{remark} \color{black}
Here, for $j \in \mathcal{J}$, we make an abuse of notation, writing $\r^\alpha$ (as well as $\R^{e,\alpha}$ and $\G^\alpha$) with $\alpha \in \{1, \ldots, d_j\}$ to refer to the variable $\r^i$ (resp. $\R^{e,i}$ and $\G^i$) where $i \in \mathcal{I}^j$ is the index of the string incident with the ``internal'' node $\N_\alpha^j$ (see also Fig. \ref{fig:notation_indices}).
\end{remark}

\begin{figure}
 \includegraphics[scale=0.6]{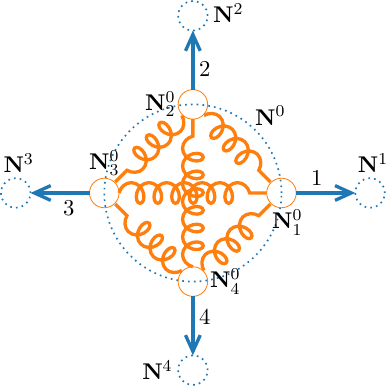}\hspace{0.5cm}
  \includegraphics[scale=0.6]{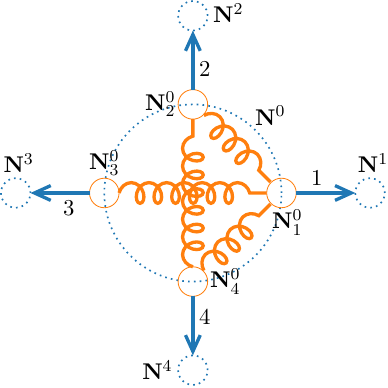}\hspace{0.5cm}
   \includegraphics[scale=0.6]{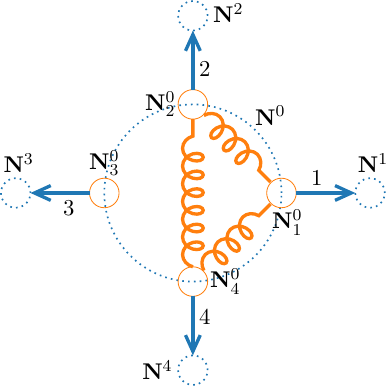}
\caption{
Diagrams of star-like networks consisting of $n=4$ strings, either with rank $n-1$ Laplacian matrix (left, middle), or with rank $< n-1$ Laplacian matrix (right).}
\label{fig:star-Laplacian}
\end{figure}

\begin{example}\label{ex-star}
We consider a star graph with $n$ strings of length $L_i = L$. We assume that the strings, at $x=0$, are connected to the only multiple node $\N^0$, the dynamics of which can be described as a local spring graph with $n$ nodes $\{\N^0_i\}_{i=1}^n$ consisting of $p$ springs ($p\geq n$).
The adjacency matrix of the spring network is given by $A=(a_{\alpha \beta})_{\alpha, \beta =1,\dots, n}$. We introduce, as in Remark \ref{vectorial}, the degree matrix $D = \mathrm{diag}\big( \sum_{\beta = 1}^n a_{\alpha \beta}, \ \alpha = 1, \ldots, n \big)$ and consequently, the discrete Laplacian ${\bf L} = D - A$. We assume that no gravitational forces are present and no other external distributed loads act on the system, such that zero is the only stretched equilibrium state. Then, the system can be rewritten simply as follows.
\begin{align}\label{star}
\begin{cases}
\r_{tt}^i(x,t)-[\G^i(\R^{e,i}_x(x)+\r^i_x(x,t))]_x = \0,
&\text{in }(0, L)\times(0,T), \, i \in \cI\\
{\bf M}\r_{tt}(0,t)={\bf K}(\r_x(0,t))-\kappa {\bf L}\r(0,t),
&t\in (0,T)\\
\r^i(L,t) = \u^i(t),
&t\in (0,T), \, i \in \cI\\
(\r^i, \r_t^i)(x,0)= (\r^{0, i},\r^{1, i})(x),
& x\in (0,L), \, i \in \cI,
\end{cases}
\end{align}
with $\r = (\r^1, \r^2, \ldots, \r^n)^T$, $\mathbf{M} = \mathrm{diag}(m_1^0, \ldots, m_n^0)$ and 
\begin{align*}
    \mathbf{K}(\r_x(0, t)) = \left(\G^1(\R^{e, 1}(0)+\r_x^1(0, t))^T, \ldots, \G^n(\R^{e, n}(0)+\r_x^n(0, t))^T \right)^T.
\end{align*}
If, in addition, the strings are all scalar in the sense that the displacements are out-of-plane only, then \eqref{star} reduces to the model discussed in \cite{WangLeugeringLi2017}.
\end{example}

\begin{remark}
{\color{black}
In Fig. \ref{fig:star-Laplacian}, we give three examples of local spring networks which all contain $n=4$ simple nodes. In fact, the rank of the Laplacian matrix is of no importance for the well-posedness of the star-like network, whereas we will see that the assumption $\mathrm{rank}(\mathbf{L}) = n-1$ allows for exact controllability of the star-like network considered in Section \ref{sec:control}. For the left and middle networks, the Laplacian matrices $\mathbf{L}_l$ and $\mathbf{L}_m$, respectively, read
\begin{align*}
{\bf L}_l=\begin{pmatrix}
3 & -1&-1& -1&\\ -1& 3& -1&-1\\
-1& -1 & 3& -1& \\
-1&-1&-1&3\\
\end{pmatrix}, \quad
{\bf L}_m=\begin{pmatrix}
2 & -1&-1& 0\\ 
-1& 2& 0&-1\\
-1& 0 & 1& 0 \\
0&-1&0&1\\
\end{pmatrix},
\end{align*}
and $\mathrm{rank}({\bf L_l}) = \mathrm{rank}({\bf L_m})=3 =n-1$. However, for the network on the right, the Laplacian matrix $\mathbf{L}_r$ reads 
\begin{align*}
{\bf L}_r=\begin{pmatrix}
1 & -1&0& 0\\ 
-1& 1& 0&0\\
0& 0 & 1& -1 \\
0&0&-1&1\\
\end{pmatrix},
\end{align*}
and one has $\mathrm{rank}({\bf L}_r)  < 3$.}
\end{remark}

Let us now pursue the discussion on the general case. 
Generally, one may put \eqref{sys-eqi} into the format
\begin{align}\label{format}
\begin{cases}
\rho_i \r^i_{tt}(x,t) = [\G^i(\R^{e,i}_x(x)+\r^i_x(x,t))]_x - \rho_i g \e,
&\text{in }(0, L_i) \times (0,T), \, i \in \cI\\
\r^{i_j}(x_{i_j,j},t) = \u^j(t),
&t\in (0,T), \, j\in \mathcal{J}^S\\
m_{i}^j\r_{tt}^i(x_{ij},t) = -[\epsilon_{ij} \mathbf{G}^i(\mathbf{R}_x^{e,i}(x_{ij}) + \r_x^i(x_{ij},t)) \\
\qquad \qquad \qquad \quad \ + \kappa_j( {\bf L}^j(\overline{\R}^{e, j} + \overline{\r}^j)(x_{\cdot, j},t))_i ],
& t \in (0,T), \, i\in\mathcal{I}^j, \, j\in \mathcal{J}^M \\
\qquad \qquad \quad \ \ =: {\bf H}^{ij}(\overline{\r}^j(x_{\cdot, j},t),\r^i_x(x_{i, j},t)) \\
(\r^i, \r_t^i)(x,0)= (\r^{0, i},\r^{1, i})(x),
&x\in (0,L), \, i\in\mathcal{I}.
\end{cases}
\end{align}
with $\overline{\r}^j$ and $\overline{\R}^{e, j}$ as in \eqref{eq:notation_rbarj}.

It is apparent that if $\r^i(x_{ij}, \cdot)\equiv \r^i_x(x_{ij}, \cdot) \equiv \0$ for all $i\in \mathcal{I}^j, \, j\in \mathcal{J}^M$, then the corresponding right-hand-side for the multiple nodes vanishes:
\begin{equation}\label{H}
{\bf H}^{ij}(\0,\0)=\0 \quad \text{for all } i\in \mathcal{I}^j, j\in\mathcal{J}^M.
\end{equation}
Moreover, we have
\begin{equation}\label{E}
[\G^i(\R^{e,i}_x+\r^i_x)]_x=\G^i_\v(\R^{e,i}_x+\r^i_x)\R^{e,i}_{xx}+\G^i_\v(\R^{e,i}_x+\r^i_x)\r^i_{xx}
\end{equation}
and, therefore, we can write
\begin{equation}\label{Eq}
\F^i(x;\r^i_x):= \rho_i g\e-\G^i_\v(\R^{e,i}_x(x)+\r^i_x)\R^{e,i}_{xx}, \quad \F^i(\cdot, \mathbf{0})=\mathbf{0}.
\end{equation}
Now, for our particular choice of the potential, we have
\begin{equation*}
\G_\v^i(\V)\v = V^i_{ss}(|\V|)\frac{\V\cd\v}{|\V|^2}\V +\frac{ V^i_s(|\V|)}{|\V|}\left[
\v-\frac{\V\cd\v}{|\V|^2}\V\right].
\end{equation*}
The matrix $\G_\v^i(\V)$, therefore, has eigenvalues $V^i_{ss}(|\V|)$
and $V^i_s(|\V|)/|\V|$ corresponding, respectively, to the eigenspaces spanned
by $\V$ and its 2-dimensional orthogonal complement.  These eigenvalues are
both positive if $|\V|>1$. An orthonormal basis of eigenvectors of $\G_\v^i(\V)$,
depending smoothly on $\V$ in the set theoretic complement of a specified
one-dimensional subspace of $\re^3$, can be chosen as follows:
$$
\frac{\V}{|\V|},\ \ \frac{M^i\V}{|M^i\V|},\ \ \frac{\V}{|\V|}\times \frac{M^i\V}{|M^i\V|},
$$
where $M^i$ is a fixed, invertible, skew-symmetric $3\times3$ matrix with the
specified subspace as nullspace and the third vector
is obtained by taking the cross product of the previous vectors.  We
restrict our attention to
{\em small perturbations of stretched equilibria} so that $|\R^{e,i}_x(x)|>1$ and
by continuity  $|\R^{e,i}_x(x)+\r^i_x(x,t)|>1$ for sufficiently small $\r^i_x$. 
Thus, the matrix
\begin{equation}\label{definite}
\G^i_\v(\R^{e,i}_x(x)+\r^i_x(x,t)) \text{ is positive definite }
\end{equation}
in a neighbourhood of the equilibrium $\R^{e,i}$. The differential equation reads in its quasilinear format as follows 
\begin{equation*}
\rho_i\r^i_{tt}-\G^i_\v(\R^{e,i}_x(x)+\r^i_x(x,t))\r^i_{xx}=\F^i(x; \r^i_x),
\end{equation*}
where $ \G^i_\v(\R^{e,i}_x(x)+\r^i_x(x,t)) $ is invertible, due to \eqref{definite}, and where \eqref{Eq} holds. This is what is needed in the further discussion. 

\begin{remark}\label{relation-to-Li}
Even though the format of the system \eqref{format} together with the definitions and properties \eqref{H}, \eqref{E}, \eqref{Eq} and \eqref{definite} is very close to the one considered in \cite{WangLeugeringLi2017}, e.g. (2.1)-(2.8), there (in \cite{WangLeugeringLi2017}) the equations are scalar, i.e. the stress depends on a scalar strain, while in the current paper the equations are vectorial, i.e. the stress, for each individual string, depends on the 3-d vector of strains. For the sake of completeness, we, therefore, provide the arguments of the proof for well-posedness in the sense of semi-global classical solutions.
\end{remark}

\medskip

\noindent \textbf{Transformation to a first order system.} We transform the system \eqref{sys-eqi}  to an equivalent  initial boundary value problem for a coupled system of first-order quasilinear
hyperbolic systems. To this end, we introduce 
\begin{align} \label{eq:wi}
\w^i = \begin{pmatrix}
\w^i_1 \\ \w^i_2 \\ \w^i_3
\end{pmatrix} 
:= \begin{pmatrix}
\r^i_x \\ \r^i_t \\ \r^i
\end{pmatrix}
\end{align}
as vector valued functions of $x$ and $t$. It is easily seen that  the equations can be rewritten as
\begin{align*}
\w^i_t + [\f^i(x,\w^i)]_x = \begin{pmatrix}
\0 \\-g\e \\ \w^{i}_2
\end{pmatrix}, \quad \f^i(x,\w^i) :=-\begin{pmatrix}
\w^i_2\\ {\rho_i}^{-1}\G^i(\R^{e,i}_x(x)+\w^i_1(x,t))\\ \0
\end{pmatrix}
\end{align*}
with $(x,t) \in [0,L_i]\times[0,T]$.
This, in turn, can be rewritten in the form of a quasilinear
hyperbolic system
\begin{equation} \label{eq:1ord_9eq}
\w^i_t + \mathbf{A}^i(x,\w^i)\w^i_x=\g^i(x,\w^i),
\end{equation}
with $\mathbf{A}^i(x,\w^i)  = \mathbf{A}^i(x,\w^i_1)$ defined by
\begin{equation*}
\mathbf{A}^i(x,\w^i) = -
\begin{pmatrix}
\0&I& \0\\\rho_i^{-1}\G_\v^i(\R^{e,i}_x(x)+\w^i_1(x,t))&\0& \0\\ \0&\0&\0
\end{pmatrix}.
\end{equation*}
and $\g^i(x,\w^i) = \g^i(x,\w^i_1,\w_2)$ defined by
\begin{equation*}
\g^i(x,\w^i) = \begin{pmatrix}\0 \\ \rho_i^{-1}\G_\v^i(\R^{e,i}_x(x) + \w^i_1)\R^{e,i}_{xx} -g\e \\ \w^{i}_2
\end{pmatrix}.
\end{equation*}
We note that $\g^i(x,\0) = \0$.
Rewriting the system \eqref{sys-eqi}, we get with 
\begin{align*}
\v^j(t) := \overline{\U}^j_t(t) +\w^{0,i_j}_2(x_{i_jj})
+t\,\rho_{i_j}^{-1} \big[\G^{i_j}(\R^{e,i_j}_x+\w^{0,i_j}_1)\big]_x\big|_{x=x_{i_jj}}\  \text{for } j\in \cJ^S,
\end{align*}
the system
\begin{equation} \label{sys-first-order} \tag{FOS}
\begin{cases}
\w^i_t + \mathbf{A}^i(x,\w^i)\w^i_x=\g^i(x,\w^i),
& \text{in } (0, L_i)\times(0, T), \, i\in\cI\\
\w^{i_j}_2(x_{i_jj},t) =\v^j(t),
&t \in (0, T),\, j\in \cJ^S\\
\eps_{ij}\G^i(\R^{e,i}_x(x_{ij})+\w^{i}_1(x_{ij},t))+m^{j}_{i}\w^{i}_{2,t}(x_{ij},t)
\\ \ + \kappa_j\Big[ \Big( \sum_{k\in\mathcal{I}^j} a^{j}_{ik}\Big)(\R^{e,i}(x_{ij})+\w^{i}_3(x_{ij},t))\\ 
 \ \ \quad \quad - \sum_{k\in\mathcal{I}^j} a^{j}_{ik}(\R^{e,k}(x_{kj})+\w^{k}_3(x_{kj},t))\Big]=0,
&t \in (0, T), \, i \in \cI^j,\, j \in \cJ^M\\
\w^i(x,0) = \w^{0,i}(x),
&x \in (0, L_i), \, i\in\cI,\\
\end{cases}
\end{equation}
with initial datum
\begin{align*}
\w^{0, i} = \begin{pmatrix}
\w^{0, i}_1 \\ \w^{0, i}_2 \\ \w^{0, i}_3
\end{pmatrix} := \begin{pmatrix}
\r^{0,i}_x \\ \r^{1,i} \\ \r^{0, i}
\end{pmatrix}.
\end{align*}
%
Introducing
\begin{align*}
\begin{aligned}
H^i(\w_1^i(x_{ij},\tau)&, \w_3^i(x_{ij},\tau),\w_3^k(x_{kj},\tau)_{k: a_{ik}^j\neq 0}):=\\ &-\frac{1}{m^{j}_{i}}
\bigg\{\eps_{ij}\G^i(\R^{e,i}_x(x_{ij})+\w^{i}_1(x_{ij},\tau))\\
& + \kappa_j \bigg[ \bigg( \sum_{k\in\mathcal{I}^j} a^{j}_{ik}\bigg)(\R^{e,i}(x_{ij})+\w^{i}_3(x_{ij},\tau)) - \sum_{k\in\mathcal{I}^j} a^{j}_{ik}(\R^{e,k}(x_{kj})+\w^{k}_3(x_{kj},\tau))\bigg] \bigg\},
\end{aligned}
\end{align*}
the system \eqref{sys-first-order} is equivalent to the non-local first order system
\begin{equation}\label{sys-first-order-nonlocal-2}\tag{NLFOS}
\begin{cases}
\w^i_t + \mathbf{A}^i(x,\w^i)\w^i_x=\g^i(x,\w^i),
&\text{in } (0,L_i)\times (0,T), \, i\in\cI\\
\w^{i_j}_2(x_{i_jj},t) =\v^j(t),
&t\in (0,T), \, j\in \cJ^S\\
\w^{i}_{2}(x_{ij},t)=\w^{i}_{2}(x_{ij},0) \\
\ +\int\limits_0^t H^i(\w_1^i(x_{ij},\tau), \w_3^i(x_{ij},\tau), \w_3^k(x_{kj},\tau)_{k: a_{ik}^j\neq 0}) d\tau,
&t\in(0, T), \, i \in \cI^j, \, j \in \cJ^M\\
\w^i(x,0) = \w^{0,i}(x), 
&x\in (0,L_i), \, i\in\cI.
\end{cases}
\end{equation}
Clearly, at the equilibrium $\R^{e,i}$, we have 
\begin{align*}
H^i(\0,\0,\0)=\0, \quad i\in \mathcal{I}_j,\, j\in\mathcal{J}^M,
\end{align*}
and our final description \eqref{sys-first-order-nonlocal-2} is now written into the format provided in \cite{WangLeugeringLi2017}.

\begin{theorem}\label{existence-semi-global}

Consider the system \eqref{sys-first-order-nonlocal-2}. Let $\R^e$ be a given stretched equilibrium of \eqref{equilibrium-nonlinear-sys}. For a specified value of $T>0$ there exist a constant $c_0>0$ (depending on $T$) such that if the initial data
\begin{align}\label{eq:ic_C1}
\w^{0,i}_1, \w^{0,i}_2, \w_3^{0,i} \in
C^1([0,L_i];\re^3)
\end{align}
and the boundary data
\begin{align*}
\v^j(t)\in C^1_0([0,T];\re^3)
\end{align*}
satisfy the $C^1$-compatibility conditions and satisfy
\begin{align}
\max
\left\{\|\w^{0,i}_1\|_1 ,\|\w^{0,i}_2\|_1, \|\w^{0,i}_3\|_1,\|\v^j\|_1\right\}_{i\in \cI, \ j\in\cJ^S} < c_0
\end{align}
there exists a unique semi-global solution
\begin{align*}
\w \in \prod_{i\in\cI}C^1([0,L_i]\times[0,T];\re^9)
\end{align*}
of \eqref{sys-first-order-nonlocal-2} depending continuously on the data. 
\end{theorem}

\begin{remark}
This theorem implies the existence and uniqueness of $C^2$ solution $\R$ of the original system \eqref{full-nonlinear-sys} in the small neighborhood of $\R^e$.
\end{remark}

\begin{proof}[Proof of Theorem \ref{existence-semi-global}]  Throughout this proof, all state vectors and coefficients depend on $(x,t)$; for brevity we sometimes write $w$ instead of $w(x,t)$, etc., when no confusion can arise.

We assume without loss of generality that  $L_i=L, \forall i\in\cI$. Otherwise one can perform a rescaling by introducing on each string a new parameter $\zeta := \fr{L}{L_i}x$ so that each string is parametrized over the same interval $[0,L]$.  We then set $\widetilde\R^i(\zeta,t) := \R^i(\fr{L_i}{L}\zeta,t)$. Then $\fr{L}{L_i} \widetilde\R^i_{\zeta}(\zeta,t) = \R^i_x(\fr{L_i}{L}\zeta,t)$ and similarly for other functions of $x$ which occur in our equations. Moreover the values $\zeta_{ij} := \fr{L}{L_i}x_{ij}$ of course now take on the values $0$ or $L$.  The system \eqref{sys-first-order-nonlocal-2} can then be rewritten in terms of the variable $\zeta \in [0,L]$ for all \textit{i}, with some added purely notational complexity. Clearly, $x_{ij}\in \{0,L\}$. Therefore, we can identify $x_{ij}$ with $L$ if $\epsilon_{ij}=1$ and with $0$ otherwise. With this initial assumption the hyperbolic system of 9 equations \eqref{eq:1ord_9eq} can now be assembled into a system in $9n$ unknowns on the domain $[0,L]\times[0,T]$
\begin{equation*}
\w_t(x,t) + \mathbf{A}\!\big(x,\w(x,t)\big)\,\w_x(x,t)=\g\!\big(x,\w(x,t)\big)
\end{equation*}
for the $n$ strings of the network, by setting
\begin{equation*}
\mathbf{A} :=\diag (\mathbf{A}^1,\mathbf{A}^2, \cdots, \mathbf{A}^n),\ \ 
\w(x,t) := \begin{pmatrix}
\w^1(x,t) \\ \w^2(x,t) \\ \vdots \\ \w^n(x,t)
\end{pmatrix}, \ \ 
\g\!\big(x,\w(x,t)\big) := \begin{pmatrix}
\g^1\!\big(x,\w^1(x,t)\big) \\ \g^2\!\big(x,\w^2(x,t)\big)\\ \vdots \\ \g^n\!\big(x,\w^n(x,t)\big)
\end{pmatrix}.
\end{equation*}
We note that $\mathbf{A}(x,\w(x,t)) = \mathbf{A}(x,\w_1(x,t))$ where $\w_1(x,t) =(\w^1_1(x,t),\w^2_1(x,t),\cdots,\w^n_1(x,t))$.
The associated initial conditions are simply
\begin{equation*}
\w(x,0) = \w^{0}(x), \quad \text{ with } \w^0 := \begin{pmatrix} \w^{0,1}\\ \w^{0,2}\\ \vdots \\ \w^{0,n} \end{pmatrix}
\in C^1([0,L];\re^{9n}),
\end{equation*}
for $x \in (0, L)$.
The conditions at simple and multiple nodes are given by
\begin{align*}
&\w^{i_j}_2(x_{i_jj},t) =\v^j(t), \quad j\in \cJ^S,\\
&\w^{i}_{2}(x_{ij},t)=\w^{i}_{2}(x_{ij},0)
+\int_0^t H^i\!\Big(\w_1^i(x_{ij},\tau), \w_3^i(x_{ij},\tau),\big(\w_3^k(x_{kj},\tau)\big)_{k: a_{ik}^j\neq 0}\Big)\, d\tau, \quad i \in \cI^j,\, j \in \cJ^M,
\end{align*}
for $t \in (0, T)$.
We shall prove that the boundary conditions both at simple nodes and multiple nodes are of the type allowing solutions to hyperbolic systems. We need to analyse the structure of the system.  Let $P^i=P^i\!\big(x,\w^i_1(x,t)\big)$ denote the matrix of the linear transformation $\rho_i^{-1}\G_\v^i\!\big(\R^{e,i}_x(x)+\w^i_1(x,t)\big)$ so that $\mathbf{A}^i$ has the block matrix structure
\begin{equation}\label{blocki} \color{black}
\mathbf{A}^i\!\big(x,\w^i_1(x,t)\big)= - \begin{pmatrix}
\0&I&\0\\ P^i\!\big(x,\w^i_1(x,t)\big)&\0&\0\\\0&\0&\0
\end{pmatrix}.
\end{equation}
We restrict our attention to
{\em small perturbations of stretched equilibria} so that $|\R^{e,i}_x(x)|>1$ and
by continuity  $|\R^{e,i}_x(x)+\w^i_1(x,t)|>1$ for sufficiently small $\w^i_1$.  Let
$Q^i\!\big(x,\w^i_1(x,t)\big)$ denote the orthogonal matrix having the vectors constructed above
with $\V=\R^{e,i}_x(x) + \w^i_1(x,t)$  as columns.   Then the symmetric matrix
$P^i\!\big(x,\w^i_1(x,t)\big)$ can be diagonalized as follows:
\begin{equation*}
P^i\!\big(x,\w^i_1(x,t)\big) = Q^i\!\big(x,\w^i_1(x,t)\big)\,[D^i\!\big(x,\w^i_1(x,t)\big)]^2\,Q^i\!\big(x,\w^i_1(x,t)\big)^T,
\end{equation*}
where $D^i\!\big(x,\w^i_1(x,t)\big) =\diag\!\big(\mu^i_1(x,\w^i_1(x,t)),\mu^i_2(x,\w^i_1(x,t)),\mu^i_3(x,\w^i_1(x,t))\big)$ is diagonal with positive diagonal entries given by
\begin{equation}\label{eigenvalues}
\begin{cases}
[\mu^i_1(x,\w^i_1(x,t))]^2= \rho^{-1}_i\,V^i_{ss}\!\Big(|\R^{e,i}_x(x) +\w^i_1(x,t)|\Big),\\[0.25em]
[\mu^i_2(x,\w^i_1(x,t))]^2 =[\mu^i_3(x,\w^i_1(x,t))]^2 =\dfrac{\rho_i^{-1}\,V^i_s\!\Big(|\R^{e,i}_x(x)+\w^i_1(x,t)|\Big)}{|\R^{e,i}_x(x)+\w^i_1(x,t)|}.
\end{cases}
\end{equation}
Let $A\in \mathbb{R}^{3p\times 3p}$, for some $p$, be given by
\begin{equation*}\color{black}
A=-
\begin{pmatrix}
\0&I&\0\\P&\0&\0\\\0&\0&\0
\end{pmatrix}.
\end{equation*}
with $P$ symmetric and positive definite, having eigenvalues $(\mu_1)^2,\dots,(\mu_p)^2$
with each $\mu_i>0$. Let $D$ denote the $p\times p$  diagonal matrix having the $\mu_i$'s
along the diagonal, and $P=QD^2Q^T$ where $Q$ is orthogonal.  Then the
matrix $A$ has eigenvalues $\pm\mu_1,\pm\mu_2,\dots,\pm\mu_p,0 $ and is
diagonalizable as follows:
\begin{equation*}
A =
S\begin{pmatrix}D&\0&\0\\\0&-D&\0\\\0&\0&\0\end{pmatrix}S^{-1},
\end{equation*}
with
\begin{align*}\color{black}
S=\begin{pmatrix}
Q&Q&\0\\
-QD&+QD&\0\\ 
\0&\0&I\end{pmatrix},\ \  
S^{-1}= \tfrac{1}{2}
\begin{pmatrix}
Q^T& -D^{-1}Q^T&\0\\
Q^T& D^{-1}Q^T&\0\\
\0&\0&I
\end{pmatrix}.
\end{align*}
Applying this to each of the diagonal blocks $\mathbf{A}^i\!\big(x,\w^i_1(x,t)\big)$ we conclude that
\begin{equation*}
\mathbf{A}\!\big(x,\w_1(x,t)\big) =
S\!\big(x,\w_1(x,t)\big)\,D\!\big(x,\w_1(x,t)\big)\,S\!\big(x,\w_1(x,t)\big)^{-1},
\end{equation*}
where
\begin{equation*}
\begin{cases}
D\!\big(x,\w_1(x,t)\big) = \diag\big(\hat D^1(x,\w^1_1(x,t)),\hat D^2(x,\w^2_1(x,t)),\dots,\hat D^n(x,\w^n_1(x,t))\big),\\
S\!\big(x,\w_1(x,t)\big) = \diag\big(S^1(x,\w^1_1(x,t)),S^2(x,\w^2_1(x,t)),\dots,S^n(x,\w^n_1(x,t))\big),\\
S\!\big(x,\w_1(x,t)\big)^{-1} = \diag\big(S^1(x,\w^1_1(x,t))^{-1},S^2(x,\w^2_1(x,t))^{-1},\dots,S^n(x,\w^n_1(x,t))^{-1}\big),\\
\end{cases}
\end{equation*}
with
\begin{equation*}
\begin{cases}
\hat D^i(x,\w^i_1(x,t)) =
\begin{pmatrix}
D^i(x,\w^i_1(x,t))&\0&\0\\
\0&-D^i(x,\w^i_1(x,t))&\0\\
\0&\0&\0
\end{pmatrix},\\
\color{black}
S^i(x,\w^i_1(x,t))=\begin{pmatrix}
Q^i(x,\w^i_1(x,t))&Q^i(x,\w^i_1(x,t))&\0\\
-\,Q^i(x,\w^i_1(x,t))D^i(x,\w^i_1(x,t))&+\,Q^i(x,\w^i_1(x,t))D^i(x,\w^i_1(x,t))&\0\\
\0&\0&I
\end{pmatrix},\\
\color{black}
S^i(x,\w^i_1(x,t))^{-1}=\tfrac{1}{2}
\begin{pmatrix}
Q^i(x,\w^i_1(x,t))^T&-D^i(x,\w^i_1(x,t))^{-1}Q^i(x,\w^i_1(x,t))^T&\0\\
Q^i(x,\w^i_1(x,t))^T&D^i(x,\w^i_1(x,t))^{-1}Q^i(x,\w^i_1(x,t))^T&\0\\
\0&\0&I
\end{pmatrix}.
\end{cases}
\end{equation*}
Premultiplying the system by $S\!\big(x,\w_1(x,t)\big)^{-1}$ we get
\begin{equation*}
S\!\big(x,\w_1(x,t)\big)^{-1}\w_t(x,t) +S\!\big(x,\w_1(x,t)\big)^{-1}
\mathbf{A}\!\big(x,\w_1(x,t)\big)\w_x(x,t)=S\!\big(x,\w_1(x,t)\big)^{-1}\g\!\big(x,\w(x,t)\big).
\end{equation*}
We denote the $k$--th row of $S\!\big(x,\w_1(x,t)\big)^{-1}$ by $\l^k\!\big(x,\w_1(x,t)\big)$.  This is a left eigenvector of $\mathbf{A}\!\big(x,\w_1(x,t)\big)$ corresponding to the $k$-th eigenvalue $\lambda_k\!\big(x,\w_1(x,t)\big)$
occurring along the diagonal of $D$. One gets the equations
\begin{equation*}
\l^k\!\big(x,\w_1(x,t)\big)\left[\w_t(x,t) +\lambda_k\!\big(x,\w_1(x,t)\big)\,\w_x(x,t)\right]=\color{black}\l^k\!\big(x,\w_1(x,t)\big)\,\g\!\big(x,\w(x,t)\big)
\end{equation*}
which is exactly of the form of the equations studied in \cite{w}. It is useful to introduce $\bxi(x,t) =S\!\big(x,\w_1(x,t)\big)^{-1}\w(x,t)$.  As in \cite{s:b} one can use the implicit function theorem to show that this change of variables is invertible near $\w=\0$. With a minor abuse of notation we also write $\bxi(x,t) = \bxi\big(x,\w(x,t)\big)$. The $k-$th component of $\bxi$ is simply $\l_k\!\big(x,\w_1(x,t)\big)\w(x,t)$, an expression commonly used in the formulation of boundary conditions.  We have
\begin{equation*}
S\!\big(x,\w_1(x,t)\big)^{-1}\w_t(x,t) +D\!\big(x,\w_1(x,t)\big)\,S\!\big(x,\w_1(x,t)\big)^{-1}\w_x(x,t)=S\!\big(x,\w_1(x,t)\big)^{-1}\g\!\big(x,\w(x,t)\big),
\end{equation*}
or, for each $i\in\cI$,
\begin{equation*}
S^i\!\big(x,\w^i_1(x,t)\big)^{-1}\w^i_t(x,t) +\hat D^i\!\big(x,\w^i_1(x,t)\big)\,S^i\!\big(x,\w^i_1(x,t)\big)^{-1}\w^i_x(x,t)=S^i\!\big(x,\w^i_1(x,t)\big)^{-1}\g^i\!\big(x,\w(x,t)\big).
\end{equation*}
Corresponding to the block structure, we write the variable $\bxi$ and modes $\bxi_+, \bxi_-$ as
\begin{align*}
\bxi(x,t) = \begin{pmatrix}
\bxi^1(x,t) \\ \bxi^2(x,t) \\ \vdots \\ \bxi^n(x,t) 
\end{pmatrix}, \quad 
\bxi_\pm(x,t) = \begin{pmatrix}
\bxi^1_\pm(x,t) \\\bxi^2_\pm(x,t) \\ \vdots\\ \bxi^n_\pm(x,t) 
\end{pmatrix},
\end{align*}
with $ \bxi^i(x,t) =S^i\!\big(x,\w^i_1(x,t)\big)^{-1}\w^i(x,t)$, and where
\begin{equation}\label{plus-minus}
\begin{split}
&\color{black} \bxi^i_+(x,t) := \fr{1}{2}\left[ Q^i\!\big(x,\w^i_1(x,t)\big)^T\,\w^i_1(x,t) - M^i\!\big(x,\w^i_1(x,t)\big)\,\w^i_2(x,t)\right],\\
& \color{black} \bxi^i_-(x,t) := \fr{1}{2}\left[ Q^i\!\big(x,\w^i_1(x,t)\big)^T\,\w^i_1(x,t) + M^i\!\big(x,\w^i_1(x,t)\big)\,\w^i_2(x,t)\right],\\
& \bxi^{i}_0(x,t):=\w^i_3(x,t)
\end{split}
\end{equation}
with $M^i\!\big(x,\w^i_1(x,t)\big):=D^i\!\big(x,\w^i_1(x,t)\big)^{-1}Q^i\!\big(x,\w^i_1(x,t)\big)^T$.  These modes correspond, respectively, to the positive and the negative eigenvalues of $\mathbf{A}^i$ and to $0$. To apply the results \textcolor{black}{of \cite{WangLeugeringLi2017}}
one needs to show that boundary conditions at $0$ and $L$ can, respectively, be rewritten \textcolor{black}{equivalently} in the form
\begin{align}
\label{eq:BC_xi+}
\bxi_+(0,t) & = \big(\mathcal{G}_0(t,\bxi_-, \bxi_0)\big)(0,t) + \int\limits_0^{t} \big(\mathcal{H}_+(t,\bxi_-,\bxi_0,\bxi_+)\big)(0,\tau)\,d\tau + \mathcal{H}_0(t),\\
\ \  \bxi_-(L,t) & = \big(\mathcal{G}_L(t, \bxi_+, \bxi_0)\big)(L,t) + \int\limits_0^{t} \big(\mathcal{H}_-(t,\bxi_-,\bxi_0,\bxi_+)\big)(L,\tau)\,d\tau +\mathcal{H}_L(t), 
\end{align}
where \textcolor{black}{$\mathcal{G}_0, \mathcal{H}_+, \mathcal{H}_0, \mathcal{G}_L, \mathcal{H}_-, \mathcal{H}_L$ are $C^1$ functions of their arguments}, $\mathcal{G}_0(t,\0) =\0$, and $\mathcal{G}_L(t,\0) =\0$.
These conditions determine the ``incoming modes'' in terms of the ``outgoing modes'' as well as boundary data. The proof is now very much like the one given in \cite{LeugeringSchmidt2011}.

We consider the situation at $0$, that at $L$ being entirely analogous. The endpoint of the $i$--th string at $0$ is either at a simple node $\N^j$ with $j\in \cJ^S$ and $i=i_j$ or it meets a multiple node $\N^j$ with $j\in\cJ^M$ and $i\in \cI^j$. We consider these two cases in the Steps 1 and 2 that follow, respectively.

\medskip

\noindent \textbf{Step 1 (simple node).} In the former case we have the boundary condition
\begin{equation}\label{eq:bc_simple}
\begin{split}
\w^i_2(0,t) = \v^j(t) = \overline{\U}^j_t(t)+
\w^{0,i}_2(0)+t\,\rho_i^{-1}\big[\G^i(\R^{e,i}_x+\w^{0,i}_1)\big]_x\big|_{x=0}
\end{split}
\end{equation}
%
At $x = 0$, we get from \eqref{plus-minus} that
\begin{align}
\bxi_+^i(0, t) + \bxi_-^i(0, t) &= Q^i(0, \w_1^i(0, t))^\intercal \w_1^i(0, t) \label{eq:changevar1}\\
-\bxi_+^i(0, t) + \bxi_-^i(0, t) &= M^i(0, \w_1^i(0, t))\w_2^i(0, t)\label{eq:changevar2}
\end{align}
Below, we use the shortened notation $C_t^1 = C^1([0, T]; \re^3)$, as well as $(C_t^1)^k$ for the product space $\Pi_{i=1}^k C_t^1$.

\medskip

\noindent \textbf{Step 1.1 (determine $\w_1^i$ using the change of variable).}
For $f \colon (C_t^1)^2 \rightarrow C_t^1$ defined by $f(\btheta, \u) = - \btheta + Q^i(0, \u)^\intercal \u$, \eqref{eq:changevar1} writes equivalently as $f\left((\bxi_+^i+\bxi_-^i)(0, \cdot), \w_1^i(0, \cdot)\right)=\0$.
Since $f(\0, \0) = \0$ and $D_\u f(\0, \0) = Q^i(0, \0)^\intercal$ is invertible, by the implicit function theorem (IFT), there exists $g \colon C_t^1 \rightarrow C_t^1$ such that, in some neightborhoods of $\0$, \eqref{eq:changevar1} is equivalent to
\begin{align*}
\w_1^i(0, \cdot) = g((\bxi_+^i+\bxi_-^i)(0, \cdot)).
\end{align*} 

\medskip

\noindent \textbf{Step 1.2 (determine $\bxi_+^i$ using the change of variable).}
Then, on some neighborhoods of $\0$, \eqref{eq:changevar2} is equivalent to 
\begin{align*}
\textcolor{black}{-\bxi_+^i(0, t) + \bxi_-^i(0, t)} = M^i(0, g((\bxi_+^i+\bxi_-^i)(0, t)))\w_2^i(0, t),
\end{align*}
which is equivalent to  $\overline{f}(\bxi_-^i(0, \cdot), \w_2^i(0, \cdot), \bxi_+^i(0, \cdot))=\0$ for $\overline{f} \colon (C_t^1)^3 \rightarrow C_t^1$ defined by $\overline{f}(\btheta_1, \btheta_2) = - \u + \btheta_1 - M^i(0, g(\u+\btheta_1))\btheta_2$. Since $\overline{f}(\0, \0, \0) = \0$ and $D_\u \overline{f}(\0, \0, \0) = \textcolor{black}{- \mathbf{I}}$ is invertible, the IFT yields that there exists $\overline{g}\colon (C_t^1)^2 \rightarrow C_t^1$ such that, on some neighborhoods of $\0$, \eqref{eq:changevar2} is equivalent to
\begin{align*}
\bxi_+^i(0, t) = \overline{g}(\bxi_-^i(0, t), \w_2^i(0, t)).
\end{align*}

\medskip

\noindent \textbf{Step 1.3 (equivalent boundary condition).} 
One may now see that the boundary condition \eqref{eq:bc_simple} is equivalent to
\begin{align*} 
\bxi_+^i(0, t) &= \overline{g}(\bxi_-^i(0, t), \v^j(t)) \\
&= \mathcal{G}_0^i(t, \bxi_-^i(0, t)) + \mathcal{H}_0^i(t),
\end{align*}
for $\mathcal{H}_0^i, \mathcal{G}_0^i$ defined by $\mathcal{H}_0^i(t) := \overline{g}(0, \v^j(t))$ and $\mathcal{G}_0^i(t, \bxi_-^i(0, t)) = \overline{g}(\bxi_-^i(0, t), \v^j(t)) - \overline{g}(0, \v^j(t))$, which indeed satisfies $\mathcal{G}_0^i(t, \0) = \0$.

\medskip 

\noindent \textbf{Step 2 (multiple node).} In the alternative case where the endpoint of the $i$-th string at $x=0$ meets a multiple node
$\N^j$ one needs to rewrite the multiple node conditions appropriately. 
These read
\begin{align} \label{eq:nodal_cond}
\begin{aligned}
&m_i^j\w_2^i(0,t) = m_i^j\w_2^{i}(0, 0) - \varepsilon_{ij} \int_0^t \mathbf{G}^i(\mathbf{R}_x^{e, i}(0)+\w_1^i(0, \tau))d\tau - q^i(t, \w_3).
\end{aligned}
\end{align}
for $q^i$ defined by
\begin{align*}
q^i(t, \w_3)=\kappa_j \int_0^t \bigg[ \bigg(\sum_{k\in \mathcal{I}^j} a_{ik}^j\bigg) (\mathbf{R}^{e,i}(0) + \w_3^i(0, \tau)) - \sum_{k\in\mathcal{I}^j}a_{ik}^j (\mathbf{R}^{e,k}(x_{kj}) - \w_3^k(x_{kj}, \tau))  \bigg]d\tau.
\end{align*}
We proceed with the following 4 steps, making repeated use of the inverse mapping theorem (IMT) and implicit function theorem:

\medskip

\noindent \textbf{Step 2.1 (invert the change of variable).} 
For $f \colon (C_t^1)^2 \rightarrow (C_t^1)^2$ defined by 
\begin{align*}
f(\x_1, \x_2) = \bigg( \tfrac{1}{2} \bigg[Q^i(0, \x_1)^\intercal \x_1 -  M^i(0, \x_1)\x_2\bigg], \tfrac{1}{2} \bigg[Q^i(0, \x_1)^\intercal \x_1 + M^i(0, \x_1)\x_2 \bigg] \bigg),
\end{align*}
the first two equations in \eqref{plus-minus} are equivalent to $\left( \bxi_+^i(0, \cdot), \bxi_-^i(0, \cdot) \right) = f(\w_1^i(0, \cdot), \w_2^i(0, \cdot))$.
One may compute that 
\begin{align*} \color{black}
(D f)(\0, \0) = \frac{1}{2} \begin{pmatrix}
 (Q_0^i)^\intercal & - (M_0^i)^\intercal \\
 (Q_0^i)^\intercal & (M_0^i)^\intercal
\end{pmatrix} \quad \text{and} \quad ((D f)(\0, \0))^{-1} = \begin{pmatrix}
Q_0^i & Q_0^i \\
-Q_0^iD_0^i & Q_0^i D_0^i
\end{pmatrix},
\end{align*}
where $Q_0^i := Q^i(0, \0)$, $M_0^i := M^i(0, \0)$ and $D_0^i := D^i(0, \0)$. Since $f(\0, \0) = (\0, \0)$ and $D f(\0, \0)$ is invertible, by the IMT, $f$ is invertible on some neighborhood of $(\0, \0)$ with inverse $g= (g_1, g_2) \colon (C_t^1)^2 \rightarrow (C_t^1)^2$, implying that \eqref{eq:changevar1}-\eqref{eq:changevar2} are equivalent to
\begin{align} 
\label{eq:inverted_changevar_w1}
\w_1^i(0, \cdot) &= g_1 \left(\bxi_+^i(0, \cdot), \bxi_-^i(0, \cdot) \right)\\
\label{eq:inverted_changevar_w2}
\w_2^i(0, \cdot) &= g_2 \left(\bxi_+^i(0, \cdot), \bxi_-^i(0, \cdot) \right).
\end{align}
Moreover,  
\begin{align*}
(D g)(\0, \0) = \begin{pmatrix}
\partial_{\y_1} g_1 & \partial_{\y_2} g_1 \\
\partial_{\y_1} g_2 & \partial_{\y_2} g_2
\end{pmatrix}(\0, \0)
= ((D f)(\0, \0))^{-1}.
\end{align*}

\medskip

\noindent \textbf{Step 2.2 (determine $\bxi_+$ using the change of variable).}
For $\overline{f}\colon (C_t^1)^2 \rightarrow (C_t^1)^2$ defined by $\overline{f}(\btheta_1, \btheta_2, \u) = g_2(\u, \btheta_2) - \btheta_1$, \eqref{eq:inverted_changevar_w2} writes equivalently as $\overline{f}(\w_2^i(0, \cdot), \bxi_-^i(0, \cdot), \bxi_+^i(0, \cdot)) = \0$.
Since $\overline{f}(\0, \0, \0) = \0$ and $D_\u \overline{f}(\0, \0) = \partial_{\y_1}g_2(\0, \0) = - Q_0^iD_0^i$ is invertible, by the IFT, there exists $\overline{g} \colon (C_t^1)^2 \rightarrow C_t^1$ such that, on some neigborhoods of $\0$,  \eqref{eq:inverted_changevar_w2} is equivalent to
\begin{align} \label{eq:extact_xi+}
\bxi_+^i(0, \cdot) = \overline{g}(\w_2^i(0, \cdot), \bxi_-^i(0, \cdot)).
\end{align}
Moreover, for $(\btheta_1, \btheta_2)$ close enough to $(\0, \0)$,
\begin{align*}
D \ \overline{g}(\btheta_1, \btheta_2) = - \left[ D_\u \ \overline{f}\left(\btheta_1, \btheta_2, \overline{g}(\btheta_1, \btheta_2)\right)\right]^{-1} D_{(\btheta_1, \btheta_2)} \ \overline{f}(\btheta_1, \btheta_2, \overline{g}(\btheta_1, \btheta_2)).
\end{align*}
Thus, as $\overline{g}(\0, \0) = \0$, and using that \textcolor{black}{$D_{\btheta_1} \overline{f} \equiv - \mathbf{I}$, $D_{\btheta_2} \overline{f}(\0, \0, \0) = D_{\y_2} g(\0, \0) = Q_0^i D_0^i$ and that $D_{\u} \overline{f}(\0, \0, \0) = D_{\y_1} g(\0, \0) = -Q_0^i D_0^i$}, we obtain
\begin{align*}
D \ \overline{g}(\0, \0) &= - \left[ D_\u \ \overline{f}\left(\0, \0, \0\right)\right]^{-1} D_{(\btheta_1, \btheta_2)} \ \overline{f}(\0, \0, \0)
= \begin{pmatrix}
-M_0^i & \mathbf{I}
\end{pmatrix}.
\end{align*}
Besides, on some neighborhoods of $\0$, \eqref{eq:inverted_changevar_w1} is equivalent to
\begin{align} \label{eq:new_eq_w1}
\w_1^i(0, \cdot) = g_1\Big(\overline{g}(\w_2^i(0, \cdot), \bxi_-^i(0, \cdot)), \bxi_-^i(0, \cdot)\Big).
\end{align}

\medskip

\noindent \textbf{Step 2.3 (determine $\w_2^i$ using the nodal condition).}
{\color{black} Then, on some neighborhoods of $\0$, \eqref{eq:changevar1}-\eqref{eq:changevar2} is equivalent to \eqref{eq:extact_xi+}-\eqref{eq:new_eq_w1}, and the nodal condition \eqref{eq:nodal_cond} is equivalent to
\begin{align}\label{eq:nodal_cond_rewritew1}
\begin{aligned}
m_i^j\w_2^i(0,t) &= m_i^j\w_2^{i}(0, 0) - q^i(t, \bxi_0)\\
&\quad - \varepsilon_{ij} \int_0^t \mathbf{G}^i\bigg(\mathbf{R}_x^{e, i}(0) + g_1\Big(\overline{g}\left(\w_2^i(0, \tau), \bxi_-^i(0, \tau)\right), \bxi_-^i(0, \tau)\Big)\bigg)d\tau.
\end{aligned}
\end{align}
Let us denote by $E^i$ and $F$ the spaces $E^i = \{\phi \in C^1([0, T]; \re^3) \colon \phi(0) = \mathbf{R}^{1, i}(0)\}$ and $F = \{\psi \in C^1([0, T]; \re^3) \colon \psi(0) = \0 \}$. For $\overline{\overline{f}} \colon (C_t^1)^{d_j+1} \times E^i \rightarrow F$ defined by
\begin{align*}
\overline{\overline{f}}(\btheta_1, \btheta_2, \u) = m_i^j (\u-\u(0)) + \varepsilon_{ij} \int_0^{\cdot} \mathbf{G}^i\bigg((\mathbf{R}_x^{e, i}(0) + g_1\Big(\overline{g}(\u, \btheta_1), \btheta_1 \Big)\bigg)d\tau + q^i(\cdot, \btheta_2),
\end{align*}
the equation \eqref{eq:nodal_cond_rewritew1} writes equivalently as $\overline{\overline{f}}(\bxi_-^i(0, \cdot), \bxi_0(x_{i,\cdot}, \cdot), \w_2^i(0, \cdot)) = \0$.
Note that $\w_2^i(0, 0) = \mathbf{r}_t^i(0, 0) = \mathbf{R}^{1, i}(0)$ due to the initial conditions.
One may compute that $D_\u \overline{\overline{f}}(\0, \0, \0) \in \mathcal{L}(E^i, F)$ has the form
\begin{align*}
\Big[\Big(&D_\u \overline{\overline{f}} (\0, \0, \0)\Big)\u_0\Big](t)\\
&= m_i^j (\u_0(t)-\u_0(0)) + \varepsilon_{ij} \int_0^{t} \mathbf{G}_\mathbf{v}^i\left(\mathbf{R}_x^{e, i}(0) + \0 \right)D_{\y_1}g_1\left(\0, \0 \right) D_{\btheta_1} \overline{g}(\0, \0)\u_0(\tau) d\tau\\
&= m_i^j (\u_0(t)-\u_0(0)) + \rho_i Q_0^i D_0^i (Q_0^i)^\intercal \int_0^t \u_0(\tau) d\tau,
\end{align*}
since $\overline{g}(\0, \0)=\0$ and $g_1(\0, \0) = \0$, and since $\varepsilon_{ij} = -1$, $\G_\v^i(\R_x^{e, i}(0)) = \rho_iP_0^i$, $D_{\y_1}g_1(\0 ,\0) = Q_0^i$ and $D_{\theta_1} \overline{g}(\0, \0) = -M_0^i$.
To show that $D_\u \overline{\overline{f}} (\0, \0, \0)$ is invertible, meaning that for any $\psi \in F$ there exists a unique $\phi \in E^i$ such that $\psi = D_\u \overline{\overline{f}} (\0, \0, \0) \phi$, let us rewrite the latter as
\begin{align} \label{eq:invert_derfbb}
\psi(t) = m_i^j (\phi(t) - \phi(0)) + \rho_i Q_0^i D_0^i (Q_0^i)^\intercal \int_0^t \phi(\tau) d\tau,
\end{align}
and denote by $Z^i$ the positive definite symmetric matrix $Z^i := \frac{\rho_i}{m_i^j} Q_0^i D_0^i (Q_0^i)^\intercal$.
Differentiating \eqref{eq:invert_derfbb} in time, and taking into account that we are looking for $\phi \in E^i$, we deduce that $\phi$ necessarily fulfills the system
\begin{align} \label{eq:invert_derfbb_DE}
\begin{cases}
\frac{\mathrm{d}}{\mathrm{d}t} \phi(t) = - Z^i \phi(t) +\frac{1}{m_i^j}\frac{\mathrm{d}}{\mathrm{d}t} \psi(t) &t \in [0, T]\\
\phi(0) =\mathbf{R}^{1, i}(0),
\end{cases}
\end{align}
whose unique solution is $\phi(t) = e^{-t Z^i} \mathbf{R}^{1, i}(0) + \frac{1}{m_i^j} \int_0^t e^{-(t-\tau) Z^i}  \frac{\mathrm{d}\psi}{\mathrm{d}t} (\tau)  d\tau$, and $D_\u \overline{\overline{f}}(\0, \0, \0)$ is thus invertible. Moreover, as $\mathbf{R}^{e,i}$ is an equilibrium,
\begin{align*}
\overline{\overline{f}}(\0, \0, \0)
&= \int_0^\cdot \bigg\{\varepsilon_{ij} \mathbf{G}^i(\mathbf{R}_x^{e,i}(0)) + \kappa_j \bigg[\bigg(\sum_{k\in \mathcal{I}^j} a_{ik}^j\bigg) \mathbf{R}^{e,i}(0) - \sum_{k\in\mathcal{I}^j}a_{ik}^j \mathbf{R}^{e,k}(x_{kj}) \bigg] \bigg\}d\tau
\end{align*}
is equal to $\0$. By the IFT, there exists $\overline{\overline{g}} \colon (C_t^1)^{d_j+1} \rightarrow E^i$ such that, on some neighborhoods, \eqref{eq:nodal_cond_rewritew1} is equivalent to 
\begin{align} \label{eq:isolate_w2}
\w_2^i(0, \cdot) = \overline{\overline{g}}(\bxi_-^i(0, \cdot), \bxi_0(0, \cdot)).
\end{align}

\medskip

\noindent \textbf{Step 2.4 (equivalent nodal condition).}
Injecting this \eqref{eq:isolate_w2} into \eqref{eq:extact_xi+}, we obtain that, on some neighborhoods of $\0$, the nodal condition \eqref{eq:nodal_cond} is equivalent to the nodal condition
\begin{align*} 
\begin{aligned}
\bxi_+^i(0, \cdot) &= \overline{g}\Big( \overline{\overline{g}}\big(\bxi_-^i(0, \cdot), \bxi_0(0, \cdot)\big), \bxi_-^i(0, \cdot) \Big)\\
&= \mathcal{G}_0^i(\bxi_0 (0, \cdot), \bxi_-^i(0, t)),
\end{aligned}
\end{align*}
for $\mathcal{G}_0^i$ defined by $\mathcal{G}_0^i(\bxi_0 (0, \cdot), \bxi_-^i(0, \cdot)) = \overline{g}\left( \overline{\overline{g}}(\bxi_-^i(0, \cdot), \bxi_0(0, \cdot)), \bxi_-^i(0, \cdot) \right)$, which indeed satisfies $\mathcal{G}_0^i(\0,\0) = \0$.}
This completes the proof of all the details needed to apply the \textcolor{black}{existence theorems presented in \cite{WangLeugeringLi2017}}
in our context and thus the proof of Theorem is complete.
\end{proof}

\section{Controllability properties of \texorpdfstring{\eqref{full-nonlinear-sys}}{systE}}
\label{sec:control}

\begin{figure}
    \centering
    \includegraphics[scale = 0.6]{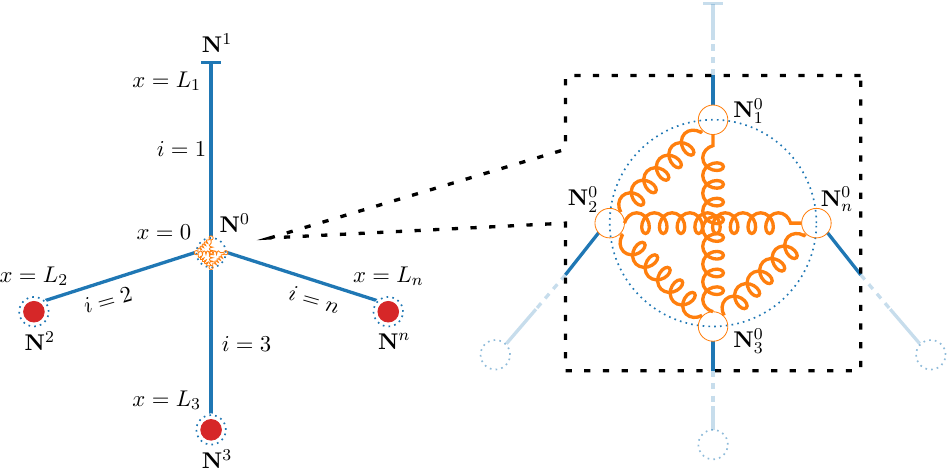}
    \caption{Star-like network clamped at the node $\mathbf{N}^1$ and controlled at the other simple nodes $\{\mathbf{N}^{i}\}_{i=2}^n$.}
    \label{fig:star_control}
\end{figure}

We consider a star-like network with $n$ strings joined together at $x=0$ at the multiple node $\N^0$, as in Example \ref{ex-star}. We exert controls at the ends for the strings labelled by $i \in \cI=\{1,...,n\}$ as follows (see also Fig. \ref{fig:star_control}). We apply homogeneous Dirichlet boundary conditions for the first string at $x = L_1$ and control the displacements of all of the other connected strings at $x = L_j$ for $j\in\cJ_c^D:=\{2,\ldots,n\}$ (Dirichlet-type controls). 
The local spring graph located at $x=0$ consists of $n$ masses $\{m_i^0\}_{i=1}^n$ connected to each other by springs of the same stiffness $\kappa$.
We recall the whole string-spring-mass system \eqref{full-nonlinear-sys} and the corresponding equilibrium system \eqref{equilibrium-nonlinear-sys} in the case of this star-like network:
 \begin{equation}\label{full-nonlinear-sys-star}
\tag{E-Star}
\begin{cases}
\rho_i \R^i_{tt}(x,t) = \G^i(\R^i_x(x,t))_x  -\rho_i g \e &\text{in }(0, L_i)\times(0,T),\, i \in \cI\nt \\
\R^{1}(L_1,t) = \0 &t\in (0,T)\nt \\
\R^{j}(L_j,t) = \U^j(t) &t\in (0,T), \,  j \in \cJ_c^D \nt \\
-\G^i(\R^i_x(0,t))+m_{i}^0\R^{i}_{tt}(0,t)\nt\\
+ \kappa\bigg[ \bigg(\sum\limits_{k\in \cI}a_{ik}\bigg)\R^{i}(0,t)-\sum\limits_{k\in \cI}a_{ik}\R^{k}(0,t)\bigg]=0 &t\in (0,T), \, i \in \cI\nt\\
(\R^i, \R^i_t)(x,0) = (\R^{0,i}, \R^{1,i})(x) &x\in (0,L_i), \, i \in \cI,\nt
\end{cases}
\end{equation}
where the stretched equilibria $\R^{e,i}(x)$ for $i\in\mathcal{I}$ satisfy
\begin{equation}\label{equilibrium-nonlinear-sys-star}
\tag{Eq-star}
\begin{cases}
\G^i(\R^{e,i}_x(x))_x  =\rho_i g \e &x \in (0, L_i), \,i \in \cI\nt \\
\R^{e,1}(L_1) = \0 \nt \\
\R^{e,j}(L_j) = \U^j &j \in \cJ_c^D \nt \\
-\G^i(\R^{e,i}_x(0))= - \kappa\bigg[ \bigg(\sum\limits_{k\in \cI}a_{ik}\bigg)\R^{e,i}(0)-\sum\limits_{k\in \cI}a_{ik}\R^{e,k}(0)\bigg] &i\in\mathcal{I}.
\end{cases}
\end{equation}
We assume the regularity and $C^2$-compatibility conditions \eqref{compatibility-1}-\eqref{compatibility-2}, where the adjacency matrix $A$, degree matrix $D$ and Laplacian matrix $\mathbf{L}$ are now taken from Example \ref{ex-star}. Furthermore, we assume that the latter satisfies $\mathrm{rank} (\mathbf{L})=n-1$, which implies that the spring graph is totally connected with at least $n-1$ springs.

\begin{remark}
The rank of $\mathbf{L}$ is in fact the difference between $n$  and the number of connected sub-graph of the string graph. Thus the graph Laplacian matrix provides a feasible means for determining the number of components in a graph, whereas there is no direct way of reading this number from the adjacency matrix.
\end{remark}

In this section, we prove Theorems \ref{local-exact-star} and \ref{global-local} below, concerning the exact controllability of the motion of a star-like network from the neighborhood of one equilibrium at time $0$ to a neighborhood of another equilibrium at time $T$. Here and below, all perturbation variables such as $\r^i(x,t)$ and their derivatives 
(e.g. $\r^i_x(x,t)$) depend on both $(x,t)$, but for brevity we sometimes omit $(x,t)$.

To make this more precise, we recall that the eigenvalues of the matrices $\mathbf{A}^i$ defined in \eqref{blocki}, are $\pm\mu_j^i(x,\w_1^i)$ for $j\in \{1,2,3\}$  given by \eqref{eigenvalues}. 
In terms of the original variables we have
\begin{equation}\label{rootsr}
\begin{cases}
\mu^i_2(x,\r^i_x) = \mu^i_3(x,\r^i_x) =
\sqrt{\frac{\rho_i^{-1}V^i_s(|\R^{e,i}_x(x)+\r^i_x|)}{|\R^{e,i}_x(x)+\r^i_x|}},\\
\mu^i_2(x,\r^i_x(x,t)) = \mu^i_3(x,\r^i_x(x,t)) = \sqrt{\frac{\rho_i^{-1}V^i_s(|\R^{e,i}_x(x)+\r^i_x(x,t)|)}{|\R^{e,i}_x(x)+\r^i_x(x,t)|}}.
\end{cases}
\end{equation}
In particular, under the assumed condition of stretched equilibria, the roots \eqref{rootsr} are strictly positive at $\w_1^i=\0,\; i\in \cI$. 
We may then define a lower bound $\overline{T}$ on the traveling time for signals entering the boundary points as Dirichlet-type controls up to their arrival at the clamped end of the first string.
Under our assumptions we can guarantee the existence of an $\epsilon_0>0$ such that 
\begin{align}
 & T>2 \overline T, \label{time1}\\
\text{and} \quad  & \overline T:= \max\left\{ T_1+\max\limits_{i=1,\dots,n}{T_i}\right\},\label{time0}
\end{align}
where each of the \emph{traveling time} for the strings labelled by $i\in \{1,\ldots, n\}$ is defined by
\begin{equation}\label{time_i}
T_i := \max_{\substack{j=1,2,3\\ x\in[0,L_i],\, t\in[0,T]}} \max_{\|\r_x^i(x,t)\|<\epsilon_0} \frac{L_i}{\mu^i_j(x,\r_x^i(x,t))}.
\end{equation}



Let the space $\H_s$ be defined by
\begin{align*}
    \H_s = \big\{ (f^i, g^i)_{i \in \mathcal{I}} \colon &(f^1, g^1) \in C^5([0, L_1]; \re^3) \times C^4 ([0, L_1]; \re^3), \\
    & \ (f^i, g^i) \in C^2([0, L_i]; \re^3) \times C^1([0, L_i]; \re^3) \text{ for }i=2, \ldots, n \, \big\}.
\end{align*}
First, in the neighborhood of any stretched equilibrium solution $(\R^e,\0)\in \H_{\text{s}}$, we show the following local result.

\begin{theorem}\label{local-exact-star}\textbf{[Local Exact Controllability]}
Let $\R^e=\{\R^{e,i}\}_{i\in\cI}$ be a stretched equilibrium solution of the system  \eqref{equilibrium-nonlinear-sys-star}. Let $T>2\overline T$, where $\overline T$ is given as \eqref{time0}. Then there exist neighborhoods $\mathcal{U}_0$ and $\mathcal{U}_1$ of $(\R^e,\0)$ in $\H_{\text{s}}$ such that for any given initial and final data 
\begin{align*}
\{(\R^{0i},\R^{1i})\} \in \cU_0 \ \text{ and }\{(\hat\R^{0i},\hat\R^{1i})\} \in \cU_1   
\end{align*}
satisfying the compatibility conditions  \eqref{compatibility-1}, \eqref{compatibility-2}. Then, one can find $(n-1)$ controls \textcolor{black}{$\U^j(t)\in C^2_0([0,T];\re^3)$ for $j \in \{2,\ldots,n\}$} such that the
corresponding solution $\R$ of \eqref{full-nonlinear-sys-star} satisfies
$$
\R^i(\cdot,T) = \hat\R^{0i}\ \text { and }  \R^i_t(\cdot,T) =
\hat\R^{1i}.
$$
\end{theorem}



Having established this local exact controllability result, we can
proceed with the following global-local exact controllability result, under the assumption that the set of equilibria is connected.

\begin{theorem}\textbf{[Global-local Exact Controllability]} \label{global-local}
Let the assumptions of Theorem~\ref{local-exact-star} hold. Given two stretched equilibrium solutions $\R^e_0(x)$ and $\R^e_1(x)$ of the system \eqref{equilibrium-nonlinear-sys-star}, there are neighborhoods $\mathcal{U}_0,\mathcal{U}_1$  of $(\R^e_0,\0), (\R^e_1,\0)$, respectively, in the state space $\H_{\text{s}}$ such that, for $T$ sufficiently large, each solution of \eqref{full-nonlinear-sys-star} in the sense of Theorem~\ref{local-exact-star} with initial conditions in $\mathcal{U}_0$ can be steered to any state in $\mathcal{U}_1$ in the given time via {\color{black} $n-1$} admissible controls.
\end{theorem}

Let us now prove these two theorems.

\subsection{Proof of the local exact controllability result}

The principal idea in boundary exact controllability of 1d-hyperbolic systems is to solve forward problems with the given initial data, backward problems with the given final data and a \textcolor{black}{sidewise Cauchy Problem from ``boundary'' to ``boundary''} (see in Fig. \ref{fig:construc}). In particular, for the latter it is convenient to interchange the spatial and time variables $x$ and $t$, and then solve Cauchy-problems from the left or the right, once the given boundary conditions have been extended to Cauchy-data there.
As we assume stretched equilibria, $\G^i_\v(\R_x^{e,i}(x))$ is
uniformly positive definite and so is
$\G_\v^i(\R_x^{e,i}(x)+\r_x^i(x,t))$ uniformly with respect to
$(x,t)$ for \textcolor{black}{$\|\r_x^i\|\le \epsilon_0$} with $\epsilon_0$ sufficiently small. The state equation for
$\r^i$ in quasilinear form can be written as
\begin{align*}
\rho_i \r_{tt}^i(x,t)=\G^i_\v(\R^{e,i}_x(x)+\r_x^i(x,t))\,\r_{xx}^i(x,t)-\rho_i g \e.
\end{align*}
Therefore,
\begin{align*}
\r_{xx}^i =\rho_i[\G^i_\v(\R_x^{e,i}+\r_x^i)]^{-1}\r_{tt}^i+\rho_ig[\G^i_\v((\R_x^{e,i}+\r_x^i)]^{-1}\e.
\end{align*}
In terms of the variables $\w^i$ introduced in \eqref{eq:wi} this reads as
\begin{align*}
\w_{1x}^i(x,t)=\rho_i[\G^i_\v(\R_x^{e,i}(x)+\w_1^i(x,t))]^{-1}\w_{2t}^i(x,t)
   +\rho_i g[\G^i_\v(\R_x^{e,i}(x)+\w_1^i(x,t))]^{-1}\e.
\end{align*}
Thus, if one has Cauchy-data $(\r^i, \r_x^i)(x,t)$ at a boundary point, for instance $x=L_i$ (i.e., $(\w_1^i,\w_2^i)(L_i,t)$ in terms of $\w^i$), and ``boundary data'' $\r^i(x,0)$ and $\r^i(x,T)$ for $x\in [0, L_i]$ (i.e., $\w_3^i(x,0)$ and $\w_3^i(x,T)$ in terms of $w^i$), then one can solve the wave equation ``from $x=L_i$ to $x=0$''. A similar statement holds for the first order system involving $\w^i$.
We emphasis that Hyperbolicity is conserved for this sidewise Cauchy Problem.

In addition, on a network, we have to focus on the interface conditions at the coupling node (here the coupling is given by the local spring-mass graph), so as to transfer all incoming information smoothly from one string to others through the coupling graph.

\begin{proof}[Proof of Theorem~\ref{local-exact-star}]
We proceed in five steps. We start with solving a time-evolution solution $\r_I = (\r_I^i)_{i\in \mathcal{I}}$ on the whole network with all given transmission conditions and given boundary conditions from $t=0$ until $t=\overline T$, which is chosen as \eqref{time0}.
The reasoning behind the choice of $\overline T$ will be made clear along the proof.

\medskip

\noindent \textbf{Step 1 (forward problem on network).}  In the first step we proceed forward from $t=0$ to \textcolor{black}{$t=\overline T$. } For each string $i\in \cI$ we define the domain $R^i_I:= [0,L_i]\times [0,\overline T]$ and for the network we set $R_I:=\prod_{i=1}^n R_I^i$. We refer to Fig. \ref{fig:construc} (left) for visualisation. \textcolor{black}{To ensure that the initial boundary value problem is well-posed, we have to provide suitable boundary data.}
The first string is fixed at $x=L_1$ with homogenous Dirichlet condition, while we impose arbitrary nonhomogeneous Dirichlet boundary conditions at $x=L_j$ for $j \in \{2,3,\ldots,n\}$. More precisely, we set $\r^j(L_j,t)=\f^j(t)$ where $\f^j(\cdot)$ are small in $C^2([0,\overline T];\re^3)$ and \textcolor{black}{satisfy the second-order compatibility conditions with the given initial condition $(\r^i,\r_t^i)(x,0)=(\phi^i,\psi^i)(x)$ where $(\phi^i,\psi^i):=(\R^{0i}-\R^{e, i}, \R^{1i})$} are themselves sufficiently small for all strings $i\in \cI$. We apply Theorem~\ref{existence-semi-global}
and obtain a unique solution $\r_I = (\r_I^i)_{i\in \mathcal{I}}$ on the forward domain $R_I$. \textcolor{black}{Furthermore, we have the asymmetric regularity $\r_I^1 \in C^5(R_I^1; \re^3)$ and $\r^i_I\in C^2(R_I^i; \re^3)$ for $i\in \cI \setminus \{1\}$.}

Then, we record  traces of state functions at some ends of each string. At the fixed simple node of the first string, along $\{x=L_1\}\times [0,\overline T]$, we take $(\ba_1,\ba_2)(t): =(\r_I^1,\r^1_{I,x})(L_1,t)$  (in fact, here $\ba_1(t)=0$ due to the homogeneous boundary condition). While at the only multiple node along $\{x=0\}\times [0,\overline T]$, we record  traces of all strings $(\bb^i_1,\bb^i_2)(t) := (\r^i_I,\r^i_{I,x})(0,t)$ for $i\in \cI$.
It is clear that 
\begin{align*}
\big(\r^i(0,t), \r^i_{tt}(0,t), \r^i_x(0,t)\big)
   = \big(\bb^i_1(t), \tfrac{\mathrm{d}^2}{\mathrm{d}t^2}\bb^i_1(t), \bb^i_2(t)\big),
   \quad t\in [0,\overline T],\; i\in\cI.
\end{align*}
satisfy the nodal conditions at the central node in \eqref{full-nonlinear-sys-star}. Moreover, all data are small in the appropriate spaces and at $x=0$ satisfy the extra regularity $\bb_1^i \in C^3([0,\overline T]; \re^3)$ (instead of only $C^2$) for all $i\in \{2, 3, \ldots, n\}$ due to the fact that $(\r_{I}^i)_{i \in \mathcal{I}}$ fulfills the interface conditions.

\begin{figure}
    \centering
    \includegraphics[scale = 0.65] {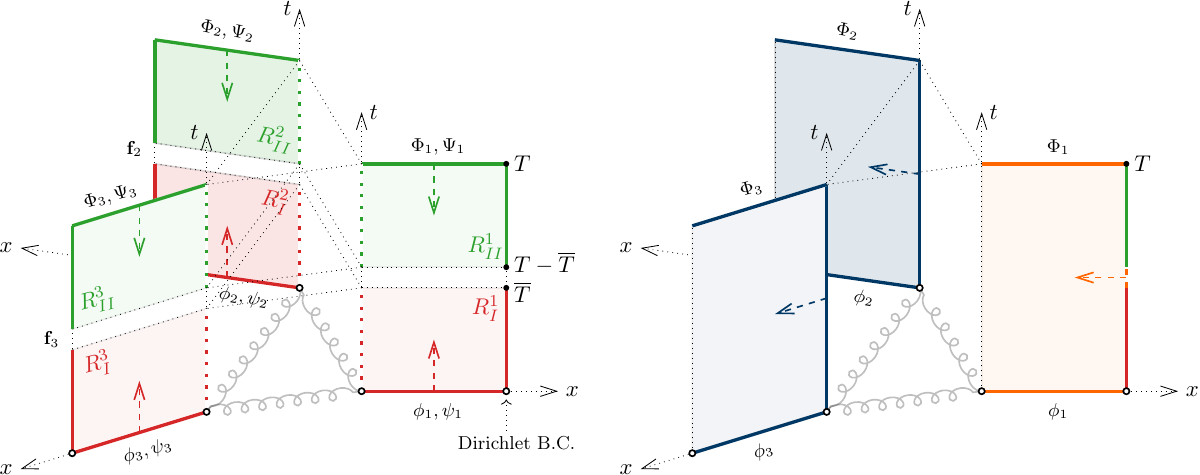}
    \caption{Constructive method with forward, backward process (left), and sidewise Cauchy problem from boundary to boundary (right)}
    \label{fig:construc}
\end{figure}

\medskip

\noindent \textbf{Step 2 (backward problem on network).} We perform the same procedure as in Step 1, but now reversing the time and progressing from the final time $T$ to $T-\overline T$, the ``final data'' being $(\r^i, \r_t^i)(x, T) = (\Phi^i, \Psi^i)(x):= (\hat{\R}^{0i} - \R^{e, i}, \hat{\R}^{1i})(x)$. More precisely, we introduce the individual domains $R_{II}^i:= [0,L_i]\times [T-\overline T, T]$ for $i \in \{1, \dots, n\}$ and for the network $R_{II}=\prod_{i=1}^n R_{II}^i$. By the same argument, a unique semi-global small solution $(\r^i_{II},\r^i_{II,x})$ of the network problem exists, and we can take traces $(\bar{\ba}_1,\bar{\ba}_2)(t) := (\r^1_{II},\r^1_{II,x})(L_1,t)$ at $\{L_1\}\times [T-\overline T, T]$ for the first string (again
$\bar{\ba}^1(t)=0$) and $(\bar{\bb}^i_1,\bar{\bb}^i_2)(t) := (\r^i_{II},\r^i_{II,x})(0,t)$
at $\{0\}\times [T-\overline T,T]$ for all strings labelled $i\in\cI$. Again, it is clear that 
\begin{align*}
\big(\r^i,\r^i_{tt},\r^i_x\big)(0,t)= \big(\bar{\bb}^i_1, \tfrac{\mathrm{d}^2}{\mathrm{d}t^2} \bar{\bb}^i_1, \bar{\bb}^i_2\big)(t),\quad t\in [T-\overline T, T], \, i \in \mathcal{I}
\end{align*}
satisfies the nodal conditions at the multiple node, and as before are small enough in appropriate spaces, and at $x=0$ we have the extra regularity $\bar{\bb}_1^i \in C^{3}([T-\overline T, T]; \re^3)$ (instead of just $C^2$) for $i \in \{2,3, \ldots, n\}$ due to the fact that $(\r_{II}^i)_{i \in \mathcal{I}}$ fulfills the interface conditions.

By the definition of $\overline T$, the domains $R_{I}$ and $R_{II}$ do not intersect. 
Following the arguments of \cite{WangLeugeringLi2017}, we connect the functions $(\ba_1, \ba_2)$ and $(\bar{\ba}_1,\bar{\ba}_2)$ (using Hermite interpolation for instance) so as to obtain functions $(\tilde{\ba}_1,\tilde{\ba}_2) \in C^5([0, T]; \re^3)\times C^4([0, T]; \re^3)$ that coincide with $(\ba_1, \ba_2)$ and $(\bar{\ba}_1,\bar{\ba}_2)$ in $[0, \overline{T}]$ and $[T-\overline{T}, T]$, respectively. This will permit us to use these Cauchy-data $(\tilde{\ba}_1,\tilde{\ba}_2)$ along $\{L_1\}\times [0,T]$ as ``initial conditions'' for the sidewise problem on the string $i=1$.

\medskip

\noindent \textbf{Step 3 (sidewise problem on string $i=1$).} We change the role of $x$ and $t$ as explained in the beginning of the proof. Let us consider the first string. The Cauchy-data just constructed is taken as data for the ``initial conditions'' $(\r^1, \r^1_x)(L_1,t)=(\tilde{\ba}^1,\tilde{\ba}^2)(t)$ for $t\in[0,T]$ at $x=L_1$, while ``boundary conditions'' at $t=0$ and $t=T$ are taken from the original initial data and given final data as $\r^1(x, 0)=\bphi^1(x)$ and $\r^1(x, T)=\bPhi^1(x)$ for $x\in [0,L_1]$. 
By \textcolor{black}{Theorem~\ref{existence-semi-global} (for a single string) to this sidewise situation}, there exist a unique solution $\r^1 \in C^5([0,L_1]\times [0,T]; \re^3)$, and we may then record the traces of $(\r^1,\r^1_x)$ along
$\{ 0\}\times [0,T]$ as follows:
\begin{align}\label{eq:solution_1}
(\tilde{\bb}^1_1, \tilde{\bb}^1_2)(t) := (\r^1,\r^1_x)(0, t),\quad t\in [0,T].
\end{align}

\smallskip

Before moving to Step 4, let us make some observations. We define the following two domains
\begin{align*}
R_{III}^1 &:= \Big\{(x, t) \colon x \in [0,L_1], t \in \Big[0, \overline T- \frac{T_1}{L_1}(L_1-x)\Big] \Big\},\\
R_{IV}^1 &:= \Big\{(x, t) \colon x \in [0,L_1], t \in \Big[T-\overline T+\frac{T_1}{L_1}(L_1-x), T\Big] \Big\}.
\end{align*}
We refer to Fig. \ref{fig:domain} for visualisation, where the time $T^*>0$ is defined as
\begin{align*}
T^* := \overline{T}-T_1.
\end{align*}
One may notice that the solution $\r_{III}^1$ to the onesided leftward problem (with ``boundary conditions'' at $t=0$) is unique in $R_{III}^1$.
Indeed, the definition of $R_{III}^1$ ensures that any characteristic curve of this problem passing by any point $(x,t)\in R_{III}^1$ is necessarily entering the domain $R_{III}^1$ at $[0, L_1]\times \{0\}$ or $ \{L_1\} \times [0, \overline{T}]$ 
Similarly, the solution $\r_{IV}^1$ to the onesided leftward problem (with ``boundary conditions'' at $t=T$) is unique in $R_{IV}^1$.
Since we chose the initial and boundary conditions for $\r_{III}^1$ and $\r_{IV}^1$ in such a way that they are also fulfilled by $\r_{I}^1$ and $\r_{II}^1$ respectively, we deduce that $\r_{I}^1 \equiv \r_{III}^1$ in $R_{III}^1$ and $\r_{II}^1 \equiv \r_{IV}^1$ in $R_{IV}^1$.
%
Therefore, at $t=0$ and $t = T$ we have
\begin{align*}
(\r_{III}^1, \r_{III, t}^1)(\cdot,0)\equiv(\bphi^1, \bpsi^1),\quad (\r_{IV}^1, \r^1_{IV, t})(\cdot,T)\equiv(\bPhi^1,\bPsi^1),\quad \text{in } [0,L_1],
\end{align*}
where $(\bphi^1, \bpsi^1)$ and $(\bPhi^1, \bPsi^1)$ are exactly the given initial and final data in Theorem \ref{local-exact-star}, while the traces of $(\r^1,\r^1_x)$ at $x=0$ satisfy
\begin{align*}
(\r_{III}^1, \r_{III, t}^1)(0,\cdot) &\equiv (\bb^1_1,\bb^1_2), \quad \text{in } [0, T^*]\\
(\r_{IV}^1, \r_{IV, t}^1)(0,\cdot) &\equiv (\bar{\bb}^1_1, \bar{\bb}^1_2), \quad \text{in }[T-T^* , T] 
\end{align*}
since one also has $\{0\}\times [0, T^*] \subset R_{III}^1$ and $\{0\} \times [T - T^*, T] \subset R_{IV}^1$.

\begin{figure}
    \centering
    \includegraphics[scale = 0.7]{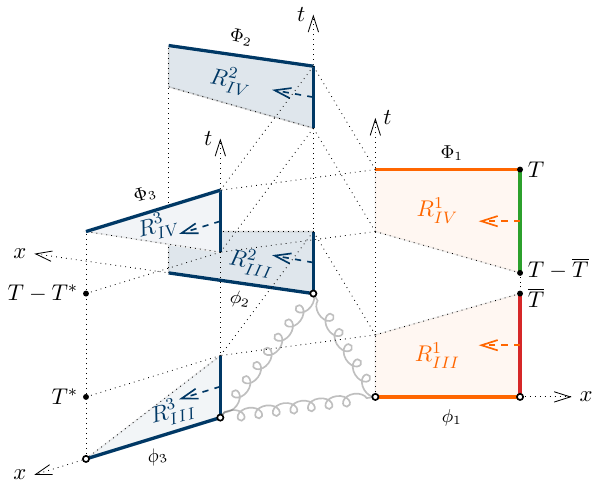}
    \caption{Different domain bounded by characteristic-lines}
    \label{fig:domain}
\end{figure}

\medskip

\noindent \textbf{Step 4 (transfer information at the junction).} In order to continue the sidewise solving process, we need to obtain the Cauchy-data $(\tilde{\bb}^i_1,\tilde{\bb}^i_2)$ of $(\r^i, \r^i_x)(0,\cdot)$ for all the \emph{remaining} strings $i\in \{2,\ldots, n\}$. The functions $\{\tilde{\bb}^i_1,\tilde{\bb}^i_2\}_{i=2}^n$ should belong to $C^2([0, T]; \re^3)\times C^1([0, T]; \re^3)$, should connect the traces of the forward and backward solutions at $x=0$:
\begin{align} \label{eq:tildbi_connect}
\begin{aligned}
&(\tilde{\bb}^i_1,\tilde{\bb}^i_2) \equiv (\bb^i_1,\bb^i_2),\quad \text{in } [0, T^*],\\
&(\tilde{\bb}^i_1,\tilde{\bb}^i_2) \equiv (\bar{\bb}^i_1, \bar{\bb}^i_2),\quad \text{in } [T-T^*, T],
\end{aligned}
\end{align}
and should satisfy the interface conditions
\begin{align}\label{eq:solve_interface}
\mathbf{M} \frac{\mathrm{d}^2}{\mathrm{d}t^2} \begin{pmatrix}
\tilde{\bb}_1^1(t)^T\\
\vdots\\
\tilde{\bb}_1^n(t)^T
\end{pmatrix} - \begin{pmatrix}
\G^1(\R^{e, 1}(0)+\tilde{\bb}_2^1(t))^T\\
\vdots\\
\G^n(\R^{e, n}(0) + \tilde{\bb}_2^n(t))^T
\end{pmatrix} + \kappa \mathbf{L} \begin{pmatrix}
\tilde{\bb}_1^1(t)^T\\
\vdots\\
\tilde{\bb}_1^n(t)^T
\end{pmatrix}=\0, \qquad t\in [0,T]
\end{align}
where $\mathbf{M}:=\diag\{m^0_1,\ldots, m^0_n\}$ and $\mathbf{L}$ are the mass and Laplacian matrices introduced in Remark \ref{vectorial}. Note that the functions $(\tilde{\bb}_1^1, \tilde{\bb}_2^1)$ are already known, given by \eqref{eq:solution_1}. Let us introduce the $n\times 3$ matrix-valued function $\Theta = (\theta_{ij})_{1\leq i \leq n, \, 1\leq j \leq 3}$ such that $\Theta(t)$ is the left-hand side of \eqref{eq:solve_interface} for all $t\in [0, T]$.

Consider the set $S := \{k\in \{2,...,n\}\colon a_{1k}\neq0\}$ containing the indexes of strings which are directly connected by a spring to the first string.
We will see below that the condition $\mathrm{rank}(\mathbf{L})=n-1$ is sufficient to ensure that $d_1 := \sum_{k \in \mathcal{I}}a_{1k}$ (which is also the cardinal of $S$) is different from zero, thereby ensuring the transfer of the information from the first string to other strings.
To find such functions $\{(\tilde{\bb}^i_1,\tilde{\bb}^i_2)\}_{i=2}^n$ solving \eqref{eq:solve_interface} with constraints \eqref{eq:tildbi_connect} (not necessarily uniquely), we proceed with two sub-steps.

\smallskip

\textbf{Step 4.1.} ~ We look into the first row of System \eqref{eq:solve_interface} which also reads $[\theta_{11}, \theta_{12}, \theta_{13}]=\0$ in $[0, T]$. Since the spring graph is connected, the first row of $\mathbf{L}$ contains at least one nonzero component other than the one on the diagonal.
Let use denote the index of one such component by $k_\circ\in S$. All three equations $[\theta_{11}, \theta_{12}, \theta_{13}] = \0$ involve $\tilde{\bb}_1^{k_\circ}$, and we may proceed the following way.
\begin{itemize}
    \item If the first string is connected to solely one other string, in other words $S = \{k_\circ\}$, then $\tilde{\bb}_1^{k_\circ}\in C^3_t$ is directly determined by the equations $[\theta_{11}, \theta_{12}, \theta_{13}] = \0$.
    
    \item Otherwise, for all indexes $k$ of strings connected to the first string other that $k_\circ$, in other words all $k \in S \setminus \{k_\circ\}$, we first find some $\tilde{\bb}_1^k \in C^3([0, T]; \re^3)$ satisfying the extension constraint \eqref{eq:tildbi_connect} (such a function is not unique). Then, the last function $\tilde{\bb}_1^{k_\circ}$, is determined by the equations $[\theta_{11}, \theta_{12}, \theta_{13}] = \0$.
\end{itemize}

\smallskip

\textbf{Step 4.2.} ~  Finally, for all strings not directly connected to the first string, which have indexes $k\in \{2, 3, \ldots, n\} \setminus S$, (similarly to Step 4.1) we first find some $\tilde{\bb}_1^k \in C^3([0, T]; \re^3)$ (not uniquely) satisfying the extension constraint \eqref{eq:tildbi_connect}. Then, the remaining $n-1$ equations of \eqref{eq:solve_interface}, which also read $[\theta_{i1}, \theta_{i2}, \theta_{i3}] = \0$ in $[0, T]$ for all $i\in \{2, 3, \ldots, n\}$ have $n-1$ unknown functions $\{\tilde{\bb}_2^i\}_{i=2}^n$. Since \eqref{definite} holds in the neighbourhood of the equilibrium $\R^{e,i}$, the functions  $\tilde{\bb}_2^i\in C^1([0, T]; \re^3)$ are uniquely determined by these equations.

\smallskip

Let us conclude this step with some remarks. The functions $\tilde{\bb}_1^{k_\circ}$ and $\{\tilde{\bb}_2^{k}\}_{k \in \{2, \ldots, n\} \setminus S}$ fulfill the constraint \eqref{eq:tildbi_connect} since they are solutions to the same system (coming from \eqref{eq:solve_interface}) as the functions ${\bb}_1^{k_\circ}$ and $\{{\bb}_2^{k}\}_{k \in \{2, \ldots, n\} \setminus S}$ (respectively) in $[0, T^*]$, and they are also solutions to the same system as the functions $\bar{\bb}_1^{k_\circ}$ and $\{\bar{\bb}_2^{k}\}_{k \in \{2, \ldots, n\} \setminus S}$ (respectively) in $[T-T^*, T]$. Both aforementioned systems having unique solutions, \eqref{eq:tildbi_connect} necessarily holds.

\smallskip

We also emphasize that, on the one hand, the regularity $\tilde{\bb}_1^{k_\circ} \in C^3_t$ is guaranteed by our asymmetric regularity assumptions on the first string (which ensures that $\tilde{\bb}_1^1 \in C^5([0, T]; \re^3)$); and on the other hand, choosing arbitrary $\tilde{\bb}_1^{k} \in C^3_t$ for $k\neq k_\circ$ is made possible by the hidden regularity enjoyed by the functions $\{\bb_1^i\}_{i=2}^n \subset C^3_t$ and $\{\bar{\bb}_1^i\}_{i=2}^n \subset C^3_t$ (observed in Steps 1 and 2, respectively).
Finally, for $k\in \{2, 3, \ldots, n\} \setminus S$, the regularity $\tilde{\bb}_2^k\in C^1([0, T]; \re^3)$ is then guaranteed by the fact that $\tilde{\bb}_1^k \in C^3([0, T]; \re^3)$ (i.e., $\frac{\mathrm{d}^2}{\mathrm{d}t^2}\tilde{\bb}_1^k \in C^1([0, T]; \re^3)$).

\medskip

\noindent \textbf{Step 5 (sidewise problem on strings $i = 2, \ldots, n$).} We now have Cauchy-data on $\{ 0\}\times [0,T]$ such that
the nodal conditions are satisfied. Therefore, we can use these as
compatible initial conditions for the strings labelled $i=2,\dots, n$
after interchanging $x$ and $t$. Thus, on the domains
$R_{IV}^i:=\{(x,t)\in [0,L_i]\times [0,T]\}$ (see in Fig.\ref{fig:domain}) we solve the initial
boundary value problems with Cauchy-data
\begin{align*}
(\r^i, \r_x^i)(0,t)=(\tilde{\bb}^i_1,\tilde{\bb}^i_2)(t),\quad t\in [0,T]
\end{align*}
nd boundary conditions
\begin{align*}
\r^i(x,0)=\bphi^i(x),\quad \r^i(x,T)=\bPhi^i(x), \quad x\in [0,L_i].
\end{align*}
By construction, the solutions are small in the sense described
above. 
Let $T^*:=\overline{T} - T_1$.
A similar uniqueness argument to that of Step 3 applies to the one-sided rightward solutions $\r_{III}^i, \r_{IV}^i$ in the respective ``lower'' and ``upper'' domains
\begin{align*}
R_{III}^i&:=\Big\{ (x,t) \colon x\in [0,L_i],\ t \in \Big[0, T^*-\frac{T_i}{L_i}x\Big] \Big\},\\
R_{IV}^i&:= \Big\{ (x,t) \colon x\in [0,L_i],\ t \in \Big[ T-T^*+\frac{T_i}{L_i}x, T\Big] \Big\}
\end{align*}
to the effect that
\begin{align*}
(\r^i, \r_t^i)(x,0) = (\bphi^i, \bpsi^i)(x), \quad (\r^i, \r_t^i)(x,T)=(\bPhi^i, \bPsi^i)(x),\quad x\in [0,L_i], i=2,...n,
\end{align*}
when
\begin{align}\label{T_star}
    T^* - T_i>0, \quad \forall i\in 2,..,n.
\end{align}

In summary, by Step 3 and Step 5, we construct a solution satisfying any given initial and final data in the neighborhoods of equilibrium in $\mathcal{U}_0$ and the system \eqref{full-nonlinear-sys-star} with the $\U^j(t)\in C^2_0([0,T];\re^3)$ for $j \in \{2,\ldots,n\}$ defined as the boundary trace of solutions. 

\end{proof}

\begin{remark}
By this constructive proof, it is clear to see the  $\overline T = T_1 + \max_{i=2,..n} \{T_i\}$ is the maximum traveling time from one control end to the clamped end on the star-like network. Furthermore, $T>2\overline T$ is optimal estimation to make sure that the forward domain in Step 1 and backward domain in Step 2 do not overlap. 
\end{remark}

\begin{remark} Here, we give more explanation of the  controllable space $\mathbf H_s$, where we require the wave on the first string (with clamped end) smoother than the other strings. 


Indeed, the boundary trace $\tilde{\bb}_1^1(t):= \mathbf r^1(0,t)$ of the solution $\mathbf r^1$ has hidden regularity at the junction $x=0$. To be specific,  $\tilde{\bb}_1^1 \in C^4([0, T]; \re^3)$, and therefore $\frac{\mathrm d}{\mathrm dt}\tilde{\bb}_1^1(t)$ and $\frac{\mathrm d^2}{\mathrm dt^2}\tilde{\bb}_1^1(t)$ are $C^2$ and $C^1$ functions of $t$, respectively. Here, we give more explanation of this hidden smooth pattern of solution due to the existence of mass (dynamical transmission conditions).

\begin{itemize}
\item \textit{Regularity of $\tilde{\bb}_1^1$ near the extremities of the time interval.} We know that the forward solution $(\r_I^i)_{i\in\mathcal{I}}$ and the backward solution $(\r_{II}^i)_{i\in\mathcal{I}}$ are both $C^2$ in their respective domains $R_I$ and $R_{II}$ and fulfill the interface conditions. Since they also coincide with $\r_{III}^1$ and $\r_{IV}^1$, respectively, in the domains $R_{III}^1$ and $R_{IV}^1$, we deduce that $\r_I^1 \in C^3(R_{III}^1; \re^3)$ and $\r_{II}^1 \in C^3(R_{IV}^1; \re^3)$. Thus, using the interface condition for $i=1$ at $x=0$, we obtain the extra-regularity $\bb_1^1 =\r^1(0,t) \in C^4([0, \overline{T}-T_1]; \re^3)$ and $\overline{\bb}_1^1=\r^1_x(0,t) \in C^4([T-\overline{T}+T_1, T]; \re^3)$, where\footnote{
An illustrative example of this argument can be given for instance when the string $i=1$ is only connected to the string $i=2$ at the junction. In this case the interface condition takes the form
\begin{align*}
m_1^0 \r_{tt}^1(0, t) = \G^1(\R^{e, 1}(0) + \r_x^1(0, t)) - \kappa \big(\R^{e, 1}(0) - \R^{e, 2}(0) + \r^1(0, t) - \r^2(0, t)\big), 
\end{align*}
for $ t \in (0, T)$. One then uses that both $\r_x^1, \r^2$ and $\r^1$ are at least $C^2$ with respect to time (the latter is in fact $C^3$) so that $\r_{tt}^1$ also has such a time regularity (rather than being only $C^1$ in time).
} we needed extra-regularity from the first string due to the presence of $\r_x^1$ in the interface condition for $i = 1$.

\item \textit{Regularity of $\tilde{\bb}_1^1$ in the middle of the time interval.} On the other hand, looking again at the leftward problem on the first string, we see that on the domain (shown as the sub-figure in the right of Fig.\ref{fig:domain})
\begin{align*}
R_{V}^1 := \Big\{ (x, t) \colon x \in [0, L_1], \ t \in \Big[ \frac{T_1}{L_1}(L_1 - x), T - T_1 - \frac{T_1}{L_1}(L_1 - x) \Big] \Big\}
\end{align*}
the solution depends only on the ``initial data'' at $x=L_1$. Assuming that this data is at least in $C^4$ with respect to $t$ (here we have homogeneous Dirichlet conditions for instance) thus ensures that $\r^1$ is also $C^4$ on $\{0\} \times [T_1, T-T_1]$, and we have $\r_x(L_1,t)$ are $C^3$ of $t$ provided the definition of $\mathbf H_s$.
\end{itemize}
\end{remark}

\subsection{Proof of the global-local exact controllability result}

We will achieve this goal in a couple of steps following essentially \cite{LeugeringSchmidt2011} with respect to the treatment of vectorial stiffness operators and \cite{WangLeugeringLi2017} with respect to the treatment of dynamic nodal conditions: we first prove the local-exact controllability result Theorem~\ref{local-exact-star}.
In the second step we use the connectedness of $\mathcal{S}_+$ and,
hence, of the set of equilibrium points. We can find a compact
$\lambda$-parametrized path $\R_\lambda^e$ of equilibrium solutions
connecting $\R^e_0$ and $\R_1^e$ and a corresponding path
$(\R^e_\lambda,\0)$ in state space. Then using a monodromy argument,
we can cover this path by sufficiently small neighborhoods such that
initial and final states in these local neighborhoods can be connected via admissible controls. In this way we can start close to the equilibrium state $(\R_0^e,\0)$ and terminate at a state close to  $(\R_1^e,\0)$ via finitely many intermediate states located in the
neighborhoods connecting $(\R_0^e,\0)$ and $(\R_1^e,\0)$. Thus we
achieve a global-local exact controllability result for the fully
nonlinear network of vibrating strings.



\begin{proof} (Theorem~\ref{global-local}.) \hspace{2mm} We now prove Theorem~\ref{global-local} for the
star-graph. As stated in the theorem, we assume two different
equilibria $\R^{e}_0, \R^e_1$ described above. These equilibria come
from different fixed boundary conditions at the simple nodes.
We assume that there is continuous  path $\R^e_\lambda$
with finite length connecting the two equilibria and hence a path $
(\R^e_\lambda,\0)$ joining the states $(\R^e_0,\0)$ and
$(\R^e_1,\0)$. Thus, given the smallness bounds needed in order to
apply Theorem~\ref{existence-semi-global}, we can cover the path with
finitely many, say $m$, such neighborhoods the centers of which are
located at an equilibrium, say $(\R^e_{\frac{k}{m}},\0)$, on the
path. The controllability times $T^k, k=1,\dots,m $ are individually
calculated according to the data in these neighborhoods according to
Theorem~\ref{local-exact-star}. Starting in the neighborhood of
$(\R^e_0,\0)$ we can reach all states $(\R^{i}(x,T^0),
\R^i_t(x,T^0)), \ i=1\dots,n$ in that neighborhood, which, in turn,
has a nonempty intersection with the next neighborhood. Therefore,
the state reached in the first step can be steered in time $T^1$ to
any state $(\R^{i}(x,T^1), \R^i_t(x,T^1)),\  i=1\dots,n$ in the
second neighborhood. This argument can be now applied $m$-times
until we arrive in the neighborhood of $(\R^e_1,\0)$. The total
control time can be estimated below by $\sum\limits_{k=1}^m T^k$.
\end{proof}

\subsection{Controllability results for the damaged model: a string-mass system coupled with less springs}

At the beginning of this section, we assume that $\text{rank}(\mathbf L) =n-1$. In the following examples, we investigate the interplay between controllability and damage (happening in form of missing springs).  Naturally, the damage problem at junction in the studied string-spring-mass system can be described in terms of the weakening of the spring stiffness $\kappa \to 0$ located at the coupling end $x = 0$ of the network, which has been investigated in \cite{LeugeringMicuRoventaWang2021} for two coupled linear strings. The control properties on this damage model might be related to the study on Laplacian matrix at the graph inside the joints. In the following cases of three mass-strings with full connection (by springs) or with reduced number of springs, we are able to see the condition $\text{rank}(\mathbf L) = n-1$ (in this case $n=3$) is \textbf{sufficient but not necessary} conditions for the controllability property for this system with $n-1$ boundary controls.



\begin{example}[The case of $\text{rank}(\mathbf L) = 2$, one missing spring] 

We consider  three strings networked by two springs, which can be used as an abstract mathematical model for flexible structures with damage at the junction.  In this toy example, there are two cases for coupling, see in Figure \ref{fig:damage}. The controllability result as  Theorem  \ref{local-exact-star} and the proof (in section 4.1) for the 'slightly' damaged model still hold by means of two controls. 

\begin{figure}[h]
    \centering
    \includegraphics[scale=0.4]{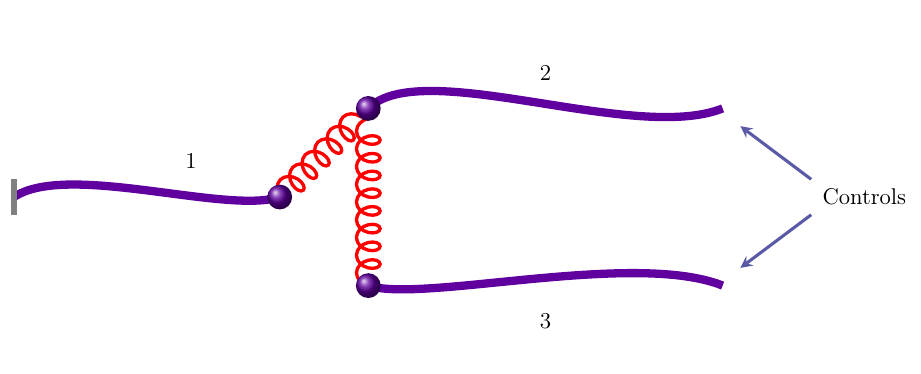}~~~ \includegraphics[scale=0.4]{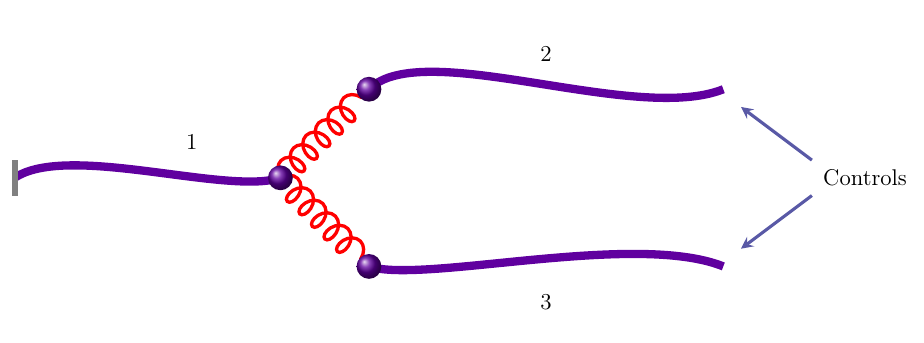} 
    \caption{three strings coupled by two springs in case 1 (left) and case 2 (right)}
    \label{fig:damage}
\end{figure}

To see this, let us have a closer look into the key step in the constructive proof, i.e. transferring information from the 1st string to the second and third one (step 4 in the proof). 

\textbf{Case 1.} The Laplacian matrix for connecting reads 
\begin{align*}
{\bf L}=\begin{pmatrix}
1 & -1 &0\\ 
-1& 2& -1\\
0& -1 & 1 \\
\end{pmatrix},
\end{align*} 
which implies the spring between string 1 and 3 is missing. Thus the transmission conditions \eqref{eq:solve_interface} with $n=3$ are written as 
\begin{subequations}
\begin{align}
&\frac{\partial^2}{\partial t^2} u^1 (0,t) = \G^1(v^1(0,t))-\kappa (u^1(0,t)-u^2(0,t)),\label{eq:damage_case1_a} \\
&\frac{\partial^2}{\partial t^2} u^2(0,t) = \G^1(v^2(0,t))-\kappa (u^2(0,t)-u^1(0,t))-\kappa(u^2(0,t)-u^3(0,t)),\label{eq:damage_case1_b}\\
&\frac{\partial^2}{\partial t^2} u^3(0,t) = \G^1(v^3(0,t))-\kappa (u^3(0,t)-u^1(0,t)),\label{eq:damage_case1_c}
\end{align}
\end{subequations}
where in order to lighten the notation, we take the end mass $m_i=1$, and denote $u^i:=\tilde{\bb^i_1} $ and $v^i:=\tilde{\bb^i_2}$ for $i=1,2,3$. The algorithm in Step 4.2 is to solve $(u^i,v^i)\in C^2_t\times C^1_t (i=2,3)$ from the above equations by given $(u^1, v^1)\in C^5_t\times C^4_t$ (which are determined by the side-wise solution on the first string in step 3). We show the steps as follows:
\begin{algorithm}[h]
\caption{(Case 1) for transferring information at junction.}
\begin{itemize}
    \item[0)]Given $(u^1(0,t), v^1(0,t))\in C^5_t[0,T]\times C^4_t[0,T]$,  thus $\frac{\partial ^2}{\partial t^2} u^1(0,t)\in C^3_t[0,T]$,
    \item[1)]Solve $u^2(0,t)\in C^3_t[0,T]$ from  \eqref{eq:damage_case1_a} when $G^1(v)$ a $C^3$ function of $v$, thus $\frac{\partial ^2}{\partial t^2} u^1(0,t)\in C^1_t[0,T]$ ,
    \item[2)] Artificially choose $u^3\in C^3_t[0,T]$ to connect the traces $u^3_f(0,t) (t\in[0,\overline T])$ and $u^3_b(0,t) (t\in [T-\overline T, T])$, which are determined by the forward and backward solutions (in step 1 and 2). Hence, $\frac{\partial ^2}{\partial t^2} u^3(0,t)\in C^1_t[0,T]$.
    \item[3)] Solve $v^2\in C^1_t[0,T]$ from \eqref{eq:damage_case1_b}, and meanwhile solve $v^3\in C^1_t[0,T]$ from \eqref{eq:damage_case1_c}. Locally, thanks to the Implicit Theorem, $(G^2)^{-1}$ and $(G^3)^{-1}$ both exist because $G^i(0)=0$, and $G^i_v(v)\neq0, v\in O(0,\epsilon)$, $i=2,3$.
\end{itemize}
\end{algorithm}

By doing this, we obtain  $(u^i,v^i)\in C^2_t[0,T]\times C^1_t[0,T] (i=2,3)$, which are taken as new Cauchy data for the side-wise solution for the second and third strings. One can easily check the side-wise solutions satisfy the whole system with given boundary and initial conditions.  
\end{example}
\bigskip

\textbf{Case 2.} In the other case, the Laplacian matrix at junction reads
\begin{align*}
{\bf L}=\begin{pmatrix}
2 & -1 &-1\\ 
-1& 1& 0\\
-1& 0 & 1 \\
\end{pmatrix},
\end{align*}
which implies the spring between 2 and 3 is missing. In details, the transmission conditions at $x=0$ are given as
\begin{subequations}
\begin{align}
&\frac{\partial^2}{\partial t^2} u^1 (0,t) = \G^1(v^1(0,t))-\kappa (u^1(0,t)-u^2(0,t)) -\kappa (u^1(0,t)-u^3(0,t)),\label{eq:damage_case2_a} \\
&\frac{\partial^2}{\partial t^2} u^2(0,t) = \G^1(v^2(0,t))-\kappa (u^2(0,t)-u^1(0,t)),\label{eq:damage_case2_b}\\
&\frac{\partial^2}{\partial t^2} u^3(0,t) = \G^1(v^3(0,t))-\kappa (u^3(0,t)-u^1(0,t)),\label{eq:damage_case2_c}
\end{align}
\end{subequations}
Similarly to case 1, for step 4, we have the following algorithm \ref{algo_2} for transferring information at junction. And the controllability is able to proved as Theorem 2.

\begin{algorithm}\label{algo_2}
\caption{(Case 2) for transferring information at junction.}
\begin{itemize}
    \item[0)]Given $(u^1(0,t), v^1(0,t))\in C^5_t[0,T]\times C^4_t[0,T]$,  thus $\frac{\partial ^2}{\partial t^2} u^1(0,t)\in C^3_t[0,T]$,
    \item[1)] Artificially choose $u^2\in C^3_t[0,T]$ to connect the traces $u^2_f(0,t) (t\in[0,\overline T])$ and $u^2_b(0,t) (t\in [T-\overline T, T])$, which are determined by the forward and backward solutions (in step 1 and 2). Hence, $\frac{\partial ^2}{\partial t^2} u^2(0,t)\in C^1_t[0,T]$.
     \item[2)] Solve $u^3(0,t)\in C^3_t[0,T]$ from  \eqref{eq:damage_case2_a} with $G^1(v)$ a $C^3$ function of $v$, thus $\frac{\partial ^2}{\partial t^2} u^3(0,t)\in C^1_t[0,T]$,
     \item[] meanwhile solve $v^2\in C^1_t[0,T]$ from \eqref{eq:damage_case2_b}
    \item[3)]Solve $v^3\in C^1_t[0,T]$ from \eqref{eq:damage_case2_c}.
\end{itemize}
\end{algorithm}

In this two cases, we only explain the step 4 in details, because we believe this is the key for the constructive method working on networks and also we observe the asymmetric controllable space in this method when the information crossing the junction. A complex smoothing pattern appears at multiple nodes when masses, springs are present. 

\begin{example}[The case of $\text{rank}(\mathbf L) = 1$, two missing spring] 
Furthermore, when the 
$\mathrm{rank}(\mathbf L)$ is reduced to 1. 
\begin{figure}[h]
    \centering
    \includegraphics[scale=0.3]{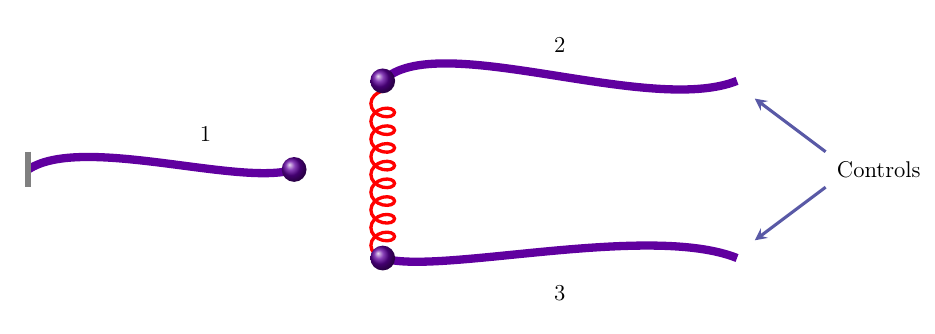} \includegraphics[scale=0.3]{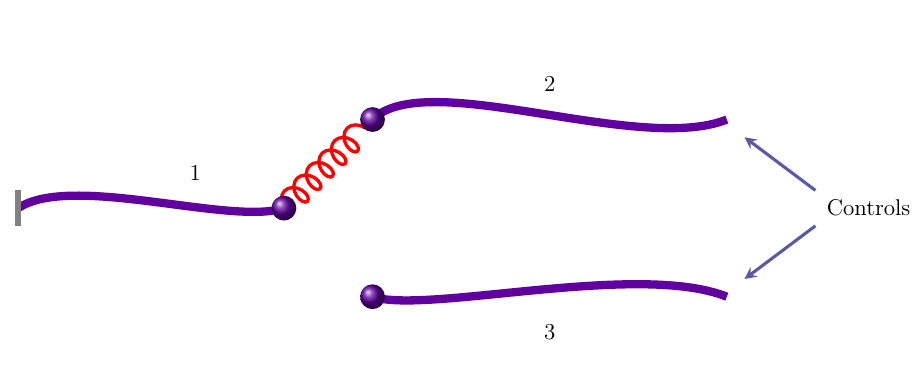} \includegraphics[scale=0.3]{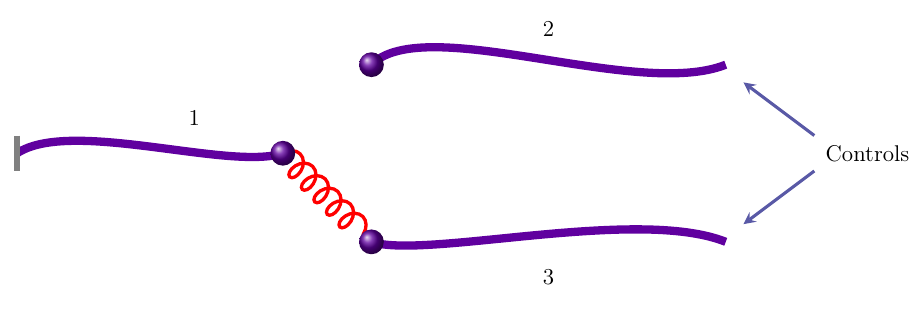} 
    \caption{Three strings coupled by one spring in case 5 (left, lack of controllability) and in case 6, 7 (middle and right).}
    \label{fig:damage_2}
\end{figure}

It is obvious that we can not always control all states on the tripod-string-networks by $2$ controls. For instance, in the case 5 (left draft in Fig. \ref{fig:damage_2}) of
\begin{align*}
{\bf L}=\begin{pmatrix}
 0 & 0 &0\\ 
 0 & 1& -1\\
 0 & -1 & 1 \\
\end{pmatrix},
\end{align*}
the clamped string (the 1st one) is not connected to others anymore, the controls given at the other two ends could never reach the clamped string, thus the system is not controllable. However, in the case 6 (left draft in Fig. \ref{fig:damage_2}) of
\begin{align*}
{\bf L}=\begin{pmatrix}
 1 & -1 &0\\ 
 -1 & 1& 0\\
 0 & 0 & 0 \\
\end{pmatrix},
\end{align*}
the system is controllable by 2 controls. The reason is that, the whole networks are divided into 2 sub-networks. The first (clamped) string and the second string are coupled, and also controllable by 1 boundary control. While the third string is free and can be controlled by the other boundary control. With this observation, we see that the rank condition, $\text{rank}= n-1$ is not necessary condition for specific 'uncoupled' networks, especially in the case that the sub-networks are equipped efficient controls separately. Hence, it is interesting to study the interplay between the controllability property of the networks and the algebra property of the Laplacian. One could investigate it in the future work. 

\end{example}




\section{Conclusion and Future Directions}
\label{sec:conclusion}

We have provided the modeling of networks of nonlinear elastic strings, where the coupling is given by networks of elastic springs and masses. Besides the genuine interest in such models, this framework can be seen as a novel modeling of potential damage in structural mechanics, where the coupling springs can be interpreted as material filaments that are subject to damage. We have discussed the well-posedness of such models and their controllability properties with particular focus on the role of the coupling spring–mass–network. The loss of controllability can now be described in terms of the controllability properties of the spring–mass–network, in addition to the well-known controllability properties of the non-linear elastic strings.\par\medskip

In this concluding section, we outline several directions for future work. First, we discuss possible extensions of the modeling framework to more general network topologies, illustrated by serial and ring-type configurations. Second, we address the inclusion of progressive damage in the coupling springs as part of the dynamical model. Finally, we reflect on methodological perspectives, in particular the use of hybrid PDE–ODE formulations, and their implications for the asymmetric controllability phenomena identified in this paper.

\subsection{Modeling Extensions}
\begin{example}[Serial connection of strings and springs]\label{ex:serial}

In the first example, we consider a serial connection of strings and springs, as represented in Figure~\ref{fig:chain_ring} (left). We take an even number of strings indexed by $i \in \{1, \dots, 2\ell\}$ for some fixed $\ell \in \{1,2,\ldots\}$, hence the network consists of $2\ell+1$ nodes. At each multiple node $\mathcal{N}^j$, $j \in \{2,3,\dots,2\ell\}$, two strings $i_1,i_2$ meet, each carrying point masses $m_{i_1}^j,m_{i_2}^j$, linked by a spring.  

The first and last nodes are simple Dirichlet-controlled boundary nodes located at $x=L$. Then the first and second strings meet at the second node at $x=0$, the second and third strings meet again at $x=L$, and so on, until the $(2\ell-1)$-th and $2\ell$-th strings meet at $x=0$. Thus, there are $2\ell-1$ multiple nodes, with $\ell$ nodes $\{2k\}_{k=1}^\ell$ at $x=0$ and $\ell-1$ nodes $\{2k-1\}_{k=2}^\ell$ at $x=L$. Each multiple node has edge degree $2$. Locally, the two incident string ends are connected in the spring graph, with adjacency, degree, and Laplacian matrices
\[
A^j=\begin{pmatrix}0&1\\1&0\end{pmatrix}, \quad
D^j=\begin{pmatrix}1&0\\0&1\end{pmatrix}, \quad
\mathbf{L}^j=\begin{pmatrix}1&-1\\-1&1\end{pmatrix}.
\]

With the notation $\overline{\mathbf{r}}^j=(\mathbf{r}^{j-1},\mathbf{r}^j)^T$, $\mathbf{M}^j$ the mass matrices, and $\mathbf{K}^j(\overline{\mathbf{r}}_x^j)$ the forces due to the strings, the nodal dynamics read
\[
\begin{aligned}
\mathbf{M}^{2k-1}\overline{\mathbf{r}}_{tt}^{2k-1}
  &= \mathbf{K}^{2k-1}(\overline{\mathbf{r}}_x^{2k-1}) - \kappa_0\mathbf{L}^{2k-1}\overline{\mathbf{r}}^{2k-1}, 
  && k=2,3,\dots,\ell,\;\; x=0, \\
\mathbf{M}^{2k}\overline{\mathbf{r}}_{tt}^{2k}
  &= -\mathbf{K}^{2k}(\overline{\mathbf{r}}_x^{2k}) - \kappa_1\mathbf{L}^{2k}\overline{\mathbf{r}}^{2k}, 
  && k=1,2,\dots,\ell,\;\; x=L,
\end{aligned}
\]
where $\kappa_0,\kappa_1$ are the spring stiffnesses at the nodes $x=0$ and $x=L$.  

This system can be rewritten as the vectorial two-point boundary value problem
\[
\begin{cases}
\mathbf{r}^i_{tt}(x,t) - [\mathbf{G}^i(\mathbf{R}^{e,i}_x(x)+\mathbf{r}^i_x(x,t))]_x = \mathbf{0}, 
& (x,t)\in(0,L)\times(0,T),\; i\in\mathcal{I}, \\
\mathbf{r}^1(L,t) = u^1(t), \;\; \mathbf{r}^{2\ell}(L,t)=u^{2\ell}(t), & t\in(0,T), \\
\mathbf{M}_0\bar{\mathbf{r}}^0_{tt}(0,t) = \mathbf{K}_0(\bar{\mathbf{r}}^0_x(0,t)) - \kappa_0\mathbf{L}_0\bar{\mathbf{r}}^0(0,t), & t\in(0,T), \\
\mathbf{M}_1\bar{\mathbf{r}}^1_{tt}(L,t) = -\mathbf{K}_1(\bar{\mathbf{r}}^1_x(L,t)) - \kappa_1\mathbf{L}_1\bar{\mathbf{r}}^1(L,t), & t\in(0,T), \\
(\mathbf{r}^i,\mathbf{r}^i_t)(x,0)=(\mathbf{r}^{0,i},\mathbf{r}^{1,i})(x), & x\in(0,L),\; i\in\mathcal{I},
\end{cases}
\]
where the block matrices $\mathbf{M}_0,\mathbf{M}_1,\mathbf{L}_0,\mathbf{L}_1$ and the force vectors $\mathbf{K}_0,\mathbf{K}_1$ are defined as in the general framework above.
\begin{figure}
\includegraphics[scale=0.75]{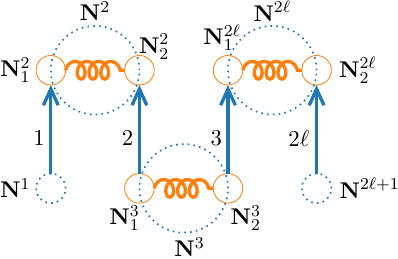} \hspace{0.75cm} \includegraphics[scale=0.75]{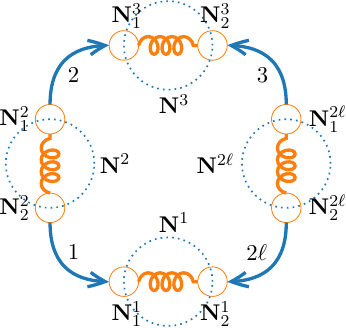}
\caption{ Diagram of the chain-like network (left), ring-like network (right) described in Examples \ref{ex:serial} and \ref{ring}, with $2\ell = 4$ strings.}
\label{fig:chain_ring}
\end{figure}

\end{example}
\begin{example}[Ring-like network]\label{ring}

We next modify Example~\ref{ex:serial} by connecting $\mathcal{N}^1$ and $\mathcal{N}^{2\\ell+1}$, so that the global network becomes a ring (cycle) with springs linking the ends of the strings; see Figure~\ref{fig:chain_ring} (right). This introduces a new multiple node at $x=0$, and no simple Dirichlet boundary nodes remain available. In this setting, controls must naturally be applied at the multiple nodes. The resulting system again takes the form of a coupled quasilinear wave equation with dynamical boundary-node conditions. In this case, we arrive at the new coupling matrix
\begin{align*}
{\bf L}_0=\begin{pmatrix}
1 & 0&0& 0&0& \dots & -1\\ 0& 1& -1&0&0& \dots & 0\\
0& -1 & 1& 0& 0&\dots & 0\\
0&0&0&1&-1&\dots &0\\0&0&0&-1&1&\dots &0\\
\ddots&\ddots&\ddots&\ddots&\ddots&\ddots&0\\
-1&0&0&0&0&\dots&1
\end{pmatrix}
\end{align*}
and new mass and strains matrices
\begin{align*}
&\mathbf{M}_0 = \mathrm{diag}(m_1^1, \mathrm{diag}(\mathbf{M}^{2k-1}, k=2, \ldots, \ell ), m_{2\ell }^1),\\
&{\mathbf{K}}_0(\r_x) = (\mathbf{G}^1(\R^{e,1}+\r^1_x), \mathbf{G}^2(\R^{e,2}+\r^2_x), \ldots, \mathbf{G}^{2\ell }(\R^{e,{2\ell }}+\r^{2\ell }_x) )^T,
\end{align*}
with the variable $\r$ defined by $\r = (\r^1, \r^2, \ldots, \r^{2\ell })^T$.
If we add controls, say at each node located at $x=0$, then we have a control vector $\u_0(t)$ such that the system reads 
\begin{align}\label{ring-eqna}
\begin{cases}
\r_{tt}^i(x,t)- [\G^i(\R^{e,i}_x(x)+\r^i_x(x,t))]_x = \0,
&\text{in }(0,L) \times (0,T), \, i \in \cI\\
{\bf M}_0\r_{tt}(0,t)={\mathbf{K}}_0(\r_x(0,t))-\kappa_0 {\bf L}_0\r(0,t)+\u_0(t),
&t\in (0,T)\\
{\bf M}_1\r_{tt}(L,t)=-{\mathbf{K}}_1(\r_x(L,t))-\kappa_1 {\bf L}_1\r(0,t),
&t\in (0,T)\\
(\r^i, \r_t^i)(x,0)= (\r^{0, i},\r^{1, i})(x),
&x\in (0, L), \, i\in\cI.
\end{cases}
\end{align}
Clearly, \eqref{ring-eqna} is again a vectorial two-point initial boundary value problem for the system of quasilinear wave equations with dynamical boundary-node-conditions, generalizing the model of \cite{WangLeugeringLi2017}. 

\begin{remark}
This example illustrates how the modeling framework extends to cycle-type networks.
Since our main controllability results concern star-like networks with Dirichlet boundary
controls at simple nodes, we do not analyze controllability for the ring configuration here.
We emphasize this model as a possible direction for future research.
\end{remark}
\end{example}

\subsection{Dynamic Damage Modeling}
In future work, we will include the dynamic evolution of damage in the coupling springs such that the spring-stiffness parameters are subject to a time evolution determined by the stresses occurring within the strings. This direction aims to capture progressive degradation in the coupling mechanism and its effect on controllability.

\subsection{Hybrid Modeling Perspective}
In this paper, we modeled the string–spring–mass network as a system of nonlinear wave equations with dynamic boundary conditions induced by point masses. This formulation is natural within the PDE framework and allowed us to apply semi-global existence theory for quasilinear hyperbolic systems. However, as pointed out during the review process, an alternative is to adopt a \emph{hybrid modeling framework}, in which the PDE dynamics of the strings and the ODE dynamics of the masses are coupled within a single evolution equation. This formulation aligns with well-known hybrid systems such as the SCOLE model, where the generator operator simultaneously encodes wave propagation in the interior and dynamical equations at the boundaries.

The hybrid viewpoint highlights that the nodal masses act as genuine internal degrees of freedom. They not only create dynamic transmission conditions but also lead to the \textbf{asymmetric controllability spaces} discovered in this work. Indeed, the smoothing effect and the resulting loss of symmetry across junctions are directly caused by the mass terms. Thus, regardless of whether the system is treated in the PDE–boundary formulation or in the hybrid PDE–ODE setting, the asymmetric exact controllability phenomenon persists.

We plan to pursue this direction in future work by developing an operator-theoretic framework that unifies the PDE and ODE dynamics into a single semigroup formulation. This will allow us to investigate stabilization and robustness properties of the asymmetric controllability phenomenon, as well as to connect with broader classes of hybrid control systems.

\section*{Acknowledgments}
The authors would like to thank Prof. Bopeng Rao (Université de Strasbourg) for helpful discussions and valuable suggestions that improved this work. This research was supported in part by the Shanghai Science and Technology Commission (STCSM) and by the Deutsche Forschungsgemeinschaft (DFG) no.~504042427.

\end{document}